\documentclass[11pt]{article}%
\usepackage[singlespacing]{setspace}
\usepackage{amsfonts}
\usepackage{graphicx}
\usepackage{amsmath}
\usepackage{amssymb}
\usepackage{layout}%
\providecommand{\U}[1]{\protect\rule{.1in}{.1in}}
\newtheorem{theorem}{Theorem}

\newtheorem{corollary}[theorem]{Corollary}

\newtheorem{example}[theorem]{Example}

\newtheorem{lemma}[theorem]{Lemma}

\newtheorem{proposition}[theorem]{Proposition}
\newtheorem{remark}[theorem]{Remark}

\newcommand\E{{\mathbb E}}
\newcommand\N{{\mathbb N}}
\newcommand\p{{\mathbb P}}
\newcommand\R{{\mathbb {R}}}
\newcommand\Z{{\mathbb {Z}}}
\newcommand\TT{{\mathbb {T}}}
\newcommand\F{{\mathcal F}}

\numberwithin{equation}{section}

\singlespace
\begin{document}

\title{On the weak invariance principle for  random fields with commuting filtrations under ${\mathbb L}^1$- projective criteria}
\author{Christophe Cuny{\thanks{Univ Brest, CNRS,  UMR 6205, Laboratoire de Math\'ematiques de Bretagne Atlantique, France, Email: christophe.cuny@univ-brest.fr}}, J\'er\^ome Dedecker{\thanks{Universit\'{e} Paris Cit\'e, Laboratoire MAP5 and CNRS UMR 8145. Email: jerome.dedecker@u-paris.fr}} and  Florence Merlev\`{e}de{\thanks{Universit\'{e} Gustave Eiffel, LAMA and CNRS
UMR 8050. Email: florence.merlevede@univ-eiffel.fr}} }\maketitle

\maketitle

\abstract{We consider a field $f \circ T_1^{i_1} \circ \cdots \circ T_d^{i_d}$ where $T_1, \dots , T_d$ are
commuting transformations, one of them  at least being ergodic. Considering the case of commuting filtrations,  we are interested by giving sufficient ${\mathbb L}^1$-projective conditions ensuring  that the 
normalized partial sums indexed by quadrants converge in distribution to a normal random variable. We also give sufficient  conditions  ensuring the weak invariance principle for the partial sums process.    For the central limit theorem (CLT), the proof  combines a truncated orthomartingale approximation with the CLT for
orthomartingales due to Voln\'y. For the functional form, a new maximal inequality is needed and is obtained via truncation techniques, blocking arguments and orthomartingale approximations.  The case of completely commuting transformations in the sense of Gordin can be handled in a similar way. Application to bounded Lipschitz functions of linear fields whose innovations have moments of a logarithmic order will be provided, as well as an application to completely commuting
endomorphisms of the $m$-torus. In the latter case, the conditions can be expressed in terms
of the ${\mathbb L}^1$-modulus of continuity of $f$.}

\medskip

\noindent{\it  MSC2020 subject classifications: }  60F05; 60F17; 60G60; 47A35; 37A05 

\section{Introduction}

\setcounter{equation}{0}

Let $(\Omega,{\mathcal A},  \p)$ be a probability space, and $T:\Omega
\mapsto \Omega$ be
 an {\it ergodic} bijective bimeasurable transformation preserving the probability $\p$. Let ${\cal F}_0 $ be a
sub-$\sigma$-algebra of ${\cal A}$ satisfying ${\mathcal F}_0
\subseteq T^{-1 }({\cal F}_0)$ and $f$ be a ${\mathbb L}^1(\p)$ real-valued centered random variable adapted to ${\cal F}_0 $.
Define then 
the stationary sequence $(f_i)_{ i \in {\mathbb Z}}$ by $f_i=f
\circ T^i$,  its associated stationary filtration $({\cal F}_i)_{ i \in {\mathbb Z}}$ by ${\mathcal F}_i={\mathcal F}_0
\circ T^{-i}$ and let $
   S_n(f) = \sum_{i=1}^{n}  f  \circ T^i  $.  
   
   Let us recall four projective criteria to get the Central Limit Theorem (CLT) for $ n^{-1/2}S_n(f)$ which are known to have a different range of applications. 
\begin{itemize}
\item The ${\mathbb L}^1$-criterion of Gordin (1973)  \cite{Go73} (see also \cite{EJ85}):
   \begin{equation} \label{L1-cond}
\sum_{i \geq 0} \E ( f_i | {\mathcal F}_0) \text{ converges in ${\mathbb L}^1(\p)$ and } \liminf_{n \rightarrow \infty} \frac{\E (|S_n(f) |)}{ \sqrt{n}} < \infty \, .
\end{equation}
\item The ${\mathbb L}^2$-criterion of Hannan (1973) \cite{Ha73}:
   \begin{equation} \label{Hannan-cond}
 f \in {\mathbb L}^2(\p) , \   \E ( f | {\mathcal F}_{- \infty})=0   \  \text{ and }  \ \sum_{i \geq 0} \Vert  \E ( f_i | {\mathcal F}_0)  -   \E ( f_i | {\mathcal F}_{-1})  \Vert_2 < \infty \, .
\end{equation}
\item The ${\mathbb L}^2$-criterion of Maxwell and Woodroofe  (2000) \cite{MW00}:
   \begin{equation} \label{mw-cond}
 f \in {\mathbb L}^2(\p) \  \text{ and }  \ \sum_{n \geq 1 }  \frac{ \Vert  \E ( S_n(f) | {\mathcal F}_0)   \Vert_2}{n^{3/2}} < \infty \, . 
\end{equation}
\item The ${\mathbb L}^1$-criterion of Dedecker and Rio (2000) \cite{DR00}:
  \begin{equation} \label{DR-cond}
 f \in {\mathbb L}^2(\p) \  \text{ and }  \  f  \E ( S_n(f) | {\mathcal F}_0)  \text{ converges in $ {\mathbb L}^1(\p)$} \, . 
\end{equation}
\end{itemize}
It has been proved by Durieu and Voln\'y \cite{DV} and Durieu \cite{Du} that these four criteria are independent.  The two ${\mathbb L}^2$-criteria \eqref{Hannan-cond} and 
\eqref{mw-cond} have been extended to the case of random fields (indexed by ${\mathbb Z}^d$, $d \geq 2$) with commuting filtrations (see Section \ref{Section1.1} for the definition of this notion).  Provided at least one of the transformations defining the random field is ergodic, the CLT  under the ${\mathbb Z}^d$-version of  Hannan's condition is given by  Voln\'y and Wang \cite{VW14} and Peligrad and Zhang \cite {PZ1} (after a previous work by Wang and Woodroofe \cite{WW13}). The extension to random fields of  Maxwell-Woodroofe's CLT has been obtained by Peligrad and Zhang \cite {PZ2}. Concerning the extension of the CLT under the  ${\mathbb L}^1$-criteria \eqref{L1-cond} or \eqref{DR-cond}, the situation is not so clear. For instance, it has been proved by Lin et al.  \cite{LMV} that in case of random fields indexed by ${\mathbb Z}^2$, the natural extension of condition \eqref{L1-cond} is not enough for the CLT to hold.  Concerning condition \eqref{DR-cond} in case of random fields, the closest result is Theorem 1 in Dedecker  \cite{D98}, but he considered 
filtrations defined with the lexicographic order (which are not commuting  in the sense given in Section \ref{Section1.1}).  

Our first objective  is then  to give an ${\mathbb L}^1$-projective criterion for the normalized partial sums associated with a stationary random field to satisfy
the central limit theorem (see our Theorem \ref{TCLCFd>2}). For random sequences (that is, when $d=1$), the approach used to prove that the CLT holds for the partial sums under conditions  \eqref{L1-cond},  \eqref{Hannan-cond} and  \eqref{mw-cond} is to prove that there exists a stationary and ergodic sequence  of martingale differences $( d_i)_{i \in {\mathbb Z}}$ in ${\mathbb L}^2$ such that $\| S_n(f) - \sum_{i=1}^n d_i\|_1 = o ( \sqrt{n})$. For the extension  to random fields of the Hannan or Maxwell and Woodroofe CLTs, the approach used is very similar and is based on an approximation of the partial sums indexed by rectangles by  orthomartingales.  However, concerning the ${\mathbb L}^1$-criteria, the paper \cite{LMV} reveals that an extension of the  ${\mathbb L}^1$-martingale $+$ coboundary decomposition as done in Gordin \cite{Go73}   does not seem  to  be the right approach.  We shall rather use an approximation by a  triangular orthomartingale in the spirit of Theorem 5.2 in \cite{HH80} (in case $d=1$), which is a different way to get an orthomartingale approximation in ${\mathbb L}^2$.  We shall also give a version of our Theorem \ref{TCLCFd>2} for completely commuting transformations in the sense of Gordin \cite{Go09}.  In Section \ref{Section2}, we shall exhibit some examples for which our ${\mathbb L}^1$-criteria  \eqref{condTCLCF1d>2} and \eqref{condTCLCF1rev} are easier to handle than the ${\mathbb L}^2$-criteria mentioned above.  Let us notice that, in case $d=1$, our ${\mathbb L}^1$-criteria  \eqref{condTCLCF1d>2} implies both criteria \eqref{L1-cond} and \eqref{DR-cond}. 

In fact, under the random fields version of the  ${\mathbb L}^2$-criteria \eqref{Hannan-cond} and 
\eqref{mw-cond}, one can prove a maximal version of the approximation of the partial sums by a stationary orthomartingale 
(see Voln\'y and Wang \cite{VW14}  under Hannan's criterion and Giraudo \cite{Gi} under Maxwell-Woodroofe's condition). Therefore the weak invariance principle for normalized partial sums indexed by rectangles follows from the weak invariance principle for orthomartingales (see Voln\'y \cite{Vo}).  Once again, the situation seems more complicated in the case of ${\mathbb L}^1$ criteria, and it seems difficult to obtain a weak invariance principle under the conditions required in our  Theorem \ref{TCLCFd>2}. In Section 4, we shall give criteria for the weak invariance principle  involving the ${\mathbb L}^1$-norm of conditional expectations of some product of random variables (see our Theorems  \ref{FCLTthm} and \ref{FCLTthmreverse}).  As usual, for the functional  CLT, the main difficulty is a suitable maximal inequality for the partial sums indexed by rectangles, which is the purpose of our Proposition \ref{Proptight}. This  new maximal inequality is  obtained via truncation techniques, blocking arguments and orthomartingale approximations.  We shall see in Subsection \ref{examplescontinued} that for the examples exhibited in Section \ref{Section2},  the weak invariance principle is satisfied under the same conditions that we found for the CLT to take place.

We shall use the following notations. For any real-valued random variable $f$ defined on $({\Omega, {\mathcal A}, {\mathbb P}})$, we denote by $Q_{f}$  the generalized inverse of the tail function $x \mapsto {\mathbb P}(|f| > x)$ and by 
$G_{f}$ the inverse of $x \mapsto \int_0^x Q_{f} (u)du$.   For any multi-integer $ {\bf i} $ of $ {\mathbb Z}^d$, we shall often use the notation $ \E_{ {\bf i} } (  \cdot)$ to mean $ \E ( \cdot | {\mathcal F}_{{\bf i}})$ for  $\sigma$-algebras indexed by $ {\bf i} $.

\section{The central limit theorem}

\setcounter{equation}{0}

\subsection{${\mathbb L}^1$ projective conditions for commuting filtrations} \label{Section1.1}

Let $d \geq 1$ and $(T_{{\bf{i}}})_{{\bf{i}} \in {\mathbb Z}^d}$ be $ {\mathbb Z}^d$ actions on  $(\Omega, {\mathcal A}, \p)$ generated by commuting invertible and measure-preserving transformations $T_{{\bf e}_q}$, $1 \leq q \leq d$. Here ${\bf e}_q$ is the vector of $ {\mathbb Z}^d$ which has $1$ at the $q$-th place and $0$ elsewhere. We assume that  at least one  of the  transformations $T_{{\bf e}_q}$, $1 \leq q \leq d$, is ergodic. By $U_{{\bf{i}}}$  we denote the operator in ${\mathbb L}^p$ ($1 \leq p \leq \infty$) defined by $ U_{{\bf{i}}} f = f \circ T_{{\bf{i}}}$, ${\bf{i}} \in {\mathbb Z}^d$.  By ${\bf{i}} \preceq {\bf{j}}$, we understand $i_k \leq j_k$ for all $1 \leq k \leq d$.  Set $S_{{\bf{n}}} (f) = \sum_{{\bf{k}} \preceq {\bf{n}}}  f \circ T_{{\bf{k}}}$ and $X_{{\bf{k}}} =   f \circ T_{{\bf{k}}}$.  For any ${\bf{n}} = (n_1, \ldots, n_d)$, set $|{\bf{n}}| = \prod_{k=1}^d n_k$. 

We suppose that  there is a $\sigma$-algebra ${\mathcal F}_{{\bf{0}}}$ such that ${\mathcal  F}_{{\bf{i}}}=T_{-{\bf{i}}} {\mathcal F}_{{\bf{0}}}$,  satisfying, for ${\bf{i}} \preceq {\bf{j}}$ ${\mathcal  F}_{{\bf{i}}} \subset {\mathcal  F}_{{\bf{j}}}$, and for an integrable $f$
\[
 \E ( \E ( f | {\mathcal  F}_{i_1, \ldots, i_d }) | {\mathcal  F}_{j_1, \ldots, j_d } )  =  \E ( f | {\mathcal  F}_{i_1 \wedge j_1, \ldots, i_d  \wedge j_d}) \, .
\]
Such a filtration  $({\mathcal  F}_{{\bf{j}}})_{{\bf{j}} \in {\mathbb Z}^d}$ is known as a   commuting filtration. 
Recall that a sequence of  orthomartingale differences $(m\circ T_{{\bf{i}}})_{\bf{i} \in {\mathbb Z}^d}$ in ${\mathbb L}^1$ adapted to a commuting filtration $ ( {\mathcal  F}_{{\bf{i}}} )_{\bf{i} \in {\mathbb Z}^d}$ means that the random variable $m$ is ${\mathcal F}_{{\bf{0}}}$ measurable and  for
each $1\leq \ell\leq d$, ${\mathbb E}(m  | {\mathcal  F}_{- \bf{e_\ell}})=0$. 

Let ${\bf e}_{i,k}=k {\bf e}_i$  the vector of $ {\mathbb Z}^d$ which has $k$ at the $i$-th place and $0$ elsewhere, and define the following weak dependence coefficients: for any $1 \leq i \leq d$ and any $k \geq 0$, 
\[
\gamma_i(k) = \Vert \E_{{\bf 0} }( X_{{\bf e}_{i,k}}) \Vert_1 \, .
\]

\begin{theorem} \label{TCLCFd>2} Assume that $f$  is in ${\mathbb L}^2$, centered and  ${\mathcal  F}_{ {\bf  0} }$-measurable, and that at least one of the transformations $T_{{\bf e}_q}$ is ergodic. Suppose in addition that 
\begin{equation} \label{condTCLCF1d>2}
   \sum_{i_1, \dots, i_d \geq 0}  \int_0^{\gamma_1(i_1) \wedge \dots \wedge \gamma_d(i_d) }  Q_f \circ G_f (u) du  < \infty \, .
\end{equation}
Then,  there exists a sequence of  orthomartingale differences $(m\circ T_{{\bf{i}}})_{\bf{i} \in {\mathbb Z}^d}$ in ${\mathbb L}^2$ such that, as  $\min_{1 \leq i \leq d}n_i \ \rightarrow \infty$, 
\begin{equation} \label{approxL2OM}
 \frac{\Vert S_{{\bf{n}}} (f)  - S_{{\bf{n}}} (m)  \Vert_2}{\sqrt{|{\bf{n}}|}}  \rightarrow 0  \, .
\end{equation}
Consequently, as  $\min_{1 \leq i \leq d}n_i \ \rightarrow \infty$, 
 \[  
\frac{1}{\sqrt{|{\bf{n}}|}} S_{{\bf{n}}} (f)   \rightarrow^{{\mathcal D}} {\mathcal N} ( 0, \Vert m \Vert_2^2)  \, .
\]
Moreover 
 \[  \Vert m \Vert_2^2 = \lim_{\min_{1 \leq i \leq d}n_i \ \rightarrow \infty}  \frac{\Vert S_{{\bf{n}}} (f)  \Vert_2^2}{|{\bf{n}}|} =  \sigma^2 (f)  \, ,
\]
where $\sigma^2 (f) = \sum_{{\bf{k}} \in {\mathbb Z}^d} {\rm Cov} ( X_{{\bf{0}}} , X_{{\bf{k}}}) $. 
\end{theorem}
\begin{remark}
If $f$ is bounded,  condition \eqref{condTCLCF1d>2} is satisfied provided
\[
  \sum_{i_1, \dots, i_d \geq 0}  \gamma_1(i_1) \wedge \dots \wedge \gamma_d(i_d)   < \infty \, .
\]
Next, if  ${\mathbb P} ( |f| >t) \leq (c/t)^r$ with $r >2$, condition \eqref{condTCLCF1d>2} holds provided
\begin{equation} \label{condsimplifiee}
  \sum_{i_1, \dots, i_d \geq 0}  \big (  \gamma_1(i_1) \wedge \dots \wedge \gamma_d(i_d)   \big )^{(r-2)/(r-1)} < \infty \, . 
\end{equation}
In addition, if   $(\gamma (k) )_{k \geq 0} \subset {\mathbb R}^+$  is a 
non-increasing sequence   such that  for any $i \in \{1, \dots, d \}$ and any $k \geq 0$,  $ \gamma_i(k) \leq \gamma (k)$, then condition \eqref{condsimplifiee} holds as soon as
\[
\sum_{k \geq 1} k^{d-1} (  \gamma ( k) )^{(r-2)/(r-1)} < \infty \, . 
\]
\end{remark}

\noindent \textbf{Proof of Theorem \ref{TCLCFd>2}.} To simplify the exposition, let us consider the case $d=2$ (the general case being identical but uses more notations). Set $T_{1,0}=T$ and $T_{0,1}=S$ and assume that the  transformations $T$ and $S$ are commuting and that $T$ (for instance) is ergodic. Let $U$ and $V$ the operators defined by 
$Uf = f \circ T$ and $Vf = f \circ S$.  We shall use the notation $\E_{a,b} (\cdot) =  \E ( \cdot | {\mathcal F}_{a,b})$ and $U^iV^j f=X_{i,j}$.  Note that, in case of $d=2$, condition \eqref{condTCLCF1d>2}
 reads as 
\begin{equation} \label{condTCLCF1}
\sum_{k, \ell \geq 0}        \int_0^{ \Vert  \E_{0,0} ( X_{k,0} )  \Vert_1 \wedge  \Vert  \E_{0,0} ( X_{0,\ell} )  \Vert_1}  Q_f \circ G_f (u) du < \infty \, .
\end{equation}
 For any $N $ a fixed positive integer, define
\[
m_N = \sum_{i =0 }^{N-1} \sum_{j=0 }^{N-1} \big (   \E_{0,0} ( U^i V^j f ) - \E_{-1,0} ( U^i V^j f )  -  \E_{0,-1} ( U^i V^j f )  +  \E_{-1,-1} ( U^i V^j f )  \big ) \, , 
\]
\[
g_{1,N} = \sum_{i =0 }^{N-1} \sum_{j=0 }^{N-1} \big (  \E_{-1,0} ( U^i V^j f )  -    \E_{-1,-1} ( U^i V^j f )  \big ) \, , \]
\[  g_{2,N} = \sum_{i =0 }^{N-1} \sum_{j=0 }^{N-1} \big (  \E_{0,-1} ( U^i V^j f )  -    \E_{-1,-1} ( U^i V^j f )  \big ) \, , 
\]
and $g_{3,N} = \sum_{i =0 }^{N-1} \sum_{j=0 }^{N-1} \E_{-1,-1} ( U^i V^j f )  $.   Hence the following decomposition is valid:
\begin{equation} \label{co-martdecwithN}
  f = m_N + (I - U)g_{1,N} + (I - V)g_{2,N}  + (I-U)(I-V)g_{3,N} + h_N \, , 
\end{equation}
where 
\begin{equation} \label{resteN}
 h_N =   \E_{0,0} ( U^N  f ) +  \E_{0,0} ( V^N  f )  -  \E_{0,0} ( U^N V^N f ) \, .
\end{equation}
For any $N$ fixed,  note that $m_N, g_{1,N}, g_{2,N}, g_{3,N}$ are in $  {\mathbb L}^2(\p)$. Moreover $(U^iV^j m_N)$ is a stationary field of orthomartingale differences, $(V^jg_{1,N})_j$ is a stationary martingale differences sequence 
with respect to the filtration $(\mathcal F_{\infty,j})_j$, and $(U^ig_{2,N})_i$ is a stationary martingale differences sequence 
with respect to the filtration $(\mathcal F_{i,\infty})_i$.    

We shall prove that 
\begin{equation} \label{approxL2OMwithN}
\lim_{N \rightarrow \infty}  \limsup_{n_1 \wedge n_2  \rightarrow \infty}\frac{\Vert S_{n_1,n_2} (f)  - S_{n_1,n_2} (m_N)  \Vert_2}{\sqrt{ n_1 n_2}}  \rightarrow 0  \, .
\end{equation}
With this aim, we first notice that 
\begin{equation*} \label{neglireste0}
\frac{1}{n_1n_2}   \big \Vert   \sum_{i=1}^{n_1}  \sum_{j=1}^{n_2}  (I - U ) U^i V^j  g_{1,N}  \big  \Vert_2^2 = 
\frac{1}{n_1}   \big \Vert   \sum_{i=1}^{n_1}   (I - U ) U^i  g_{1,N}  \big  \Vert_2^2 \leq  \frac{4}{n_1}  \Vert  g_{1,N}  \big  \Vert^2_2 \, , 
\end{equation*}
which is going to zero as $n_1 \rightarrow \infty$. Similarly
\begin{equation*} \label{neglireste00}
\frac{1}{n_1n_2}   \big \Vert   \sum_{i=1}^{n_1}  \sum_{j=1}^{n_2}  (I - V ) U^i V^j  g_{2,N}  \big  \Vert_2^2  \leq  \frac{4}{n_2}  \Vert  g_{2,N}  \big  \Vert^2_2 \, , 
\end{equation*}
which is going to zero as $n_2 \rightarrow \infty$.  In addition
\begin{equation*} \label{neglireste000}
\frac{1}{n_1n_2}   \big \Vert   \sum_{i=1}^{n_1}  \sum_{j=1}^{n_2}  (I-U)  (I - V ) U^i V^j  g_{3,N}  \big  \Vert_2^2  \leq  \frac{16}{n_1 n_2}  \Vert  g_{3,N}  \big  \Vert^2_2 \, , 
\end{equation*}
 Hence, \eqref{approxL2OMwithN} will follow if one  can prove that 
\begin{equation} \label{neglireste}
\lim_{N \rightarrow \infty} \limsup_{n_1 \wedge n_2 \rightarrow \infty}\frac{1}{n_1n_2}  \Vert  S_{n_1,n_2} (h_N)  \Vert_2^2=0 \, . 
\end{equation}
Let $h_{N,2} =\E_{0,0} ( V^N  f )$.  Note first that 
\begin{multline*}
 \Vert  S_{n_1,n_2} (h_{N,2})  \Vert_2^2  = \sum_{k_1,k_2=1}^{n_1} \sum_{\ell_1,\ell_2=1}^{n_2}  \E \big ( \E_{k_1,\ell_1} ( X_{k_1,\ell_1 +N} )  \E_{k_2,\ell_2} ( X_{k_2,\ell_2 +N} )  \big ) \\
 \leq 2     \sum_{k_1=1}^{n_1}  \sum_{k_2=k_1}^{n_1}    \sum_{\ell_1=1}^{n_2}  \sum_{\ell_2=\ell_1}^{n_2}  \Big |  \E \big ( \E_{k_1,\ell_1} ( X_{k_1,\ell_1 +N} )  \E_{k_2,\ell_2} ( X_{k_2,\ell_2 +N} )  \big ) \Big |  \\
 + 2  \sum_{k_1=1}^{n_1}  \sum_{k_2=k_1}^{n_1}    \sum_{\ell_1=1}^{n_2}  \sum_{\ell_2=1}^{\ell_1}  \Big |  \E \big ( \E_{k_1,\ell_1} ( X_{k_1,\ell_1 +N} )  \E_{k_2,\ell_2} ( X_{k_2,\ell_2 +N} )  \big ) \Big | \, .
\end{multline*}
By the properties of the conditional expectation, if $k_1 \leq k_2$ and $\ell_1 \leq \ell_2$,
\[
 \E \big ( \E_{k_1,\ell_1} ( X_{k_1,\ell_1 +N} )  \E_{k_2,\ell_2} ( X_{k_2,\ell_2 +N} )  \big ) =  \E \big ( \E_{k_1,\ell_1} ( X_{k_1,\ell_1 +N} ) \E_{k_1,\ell_1} ( X_{k_2,\ell_2 +N} )  \big ) \, .
\]
Next, again by the properties of the conditional expectation and the fact that the filtrations are completely commuting, we derive that for $k_1 \leq k_2$ and $\ell_1 \geq \ell_2$, 
\begin{multline*}
\E \big ( \E_{k_1,\ell_1} ( X_{k_1,\ell_1 +N} )  \E_{k_2,\ell_2} ( X_{k_2,\ell_2 +N} )  \big )  = \E \big ( \E_{k_1,\ell_2} ( X_{k_1,\ell_1 +N} )  \E_{k_2,\ell_2} ( X_{k_2,\ell_2 +N} )  \big ) \\
= \E \big ( \E_{k_1,\ell_2} ( X_{k_1,\ell_1 +N} )  \E_{k_1,\ell_2} ( X_{k_2,\ell_2 +N} )  \big ) \, .
\end{multline*}
So, overall for any  positive integers $k_1,k_2, \ell_1, \ell_2$,
\begin{equation} \label{controlcondexp}
\E \big ( \E_{k_1,\ell_1} ( X_{k_1,\ell_1 +N} )  \E_{k_2,\ell_2} ( X_{k_2,\ell_2 +N} )  \big ) 
= \E \big ( \E_{k_1\wedge k_2, \ell_1 \wedge \ell_2} ( X_{k_1 ,\ell_1 +N} ) \E_{k_1 \wedge k_2, \ell_1 \wedge \ell_2} ( X_{k_2,\ell_2 +N} )  \big ) \, .
\end{equation}
Hence, by using \eqref{controlcondexp} and by stationarity, we derive 
\begin{multline*}
 \Vert  S_{n_1,n_2} (h_{N,2})  \Vert_2^2 
 \leq 2  n_1 n_2  \sum_{k=0}^{n_1-1}  \sum_{\ell=0}^{n_2-1}   \big |  \E \big ( \E_{0,0} ( X_{0,N} )  \E_{0,0} ( X_{k,\ell +N} )  \big ) \big |  \\
 + 2  n_1 n_2 \sum_{k=0}^{n_1-1}     \sum_{\ell=0}^{n_2-1}   \big |  \E \big ( \E_{0,0} ( X_{0,\ell +N} )  \E_{0,0} ( X_{k,N} )  \big ) \big | \, .
\end{multline*}
Therefore
\begin{multline*}
 \limsup_{n_1 \wedge n_2 \rightarrow \infty}\frac{1}{n_1n_2}  \Vert  S_{n_1,n_2} (h_{N,2})  \Vert_2^2  \\ \leq 
2   \sum_{k \geq 0}     \sum_{\ell \geq 0}  \Big \{   \big |  \E \big ( \E_{0,0} ( X_{0,N} )  \E_{0,0} ( X_{k,\ell +N} )  \big ) \big | +   \big |  \E \big ( \E_{0,0} ( X_{0,\ell +N} )  \E_{0,0} ( X_{k,N} )  \big ) \big |  \Big \}  \, .
\end{multline*}
By Proposition 1 in \cite{DD03} and stationarity, we have
\begin{multline*}
\big \vert  \E \big ( \E_{0,0} ( X_{0,N} )  \E_{0,0} ( X_{k,\ell +N} )  \big )  \big  \vert 
\leq \int_0^{\Vert  \E_{0,0} ( X_{k,\ell} ) \Vert_1} Q_{ |\E_{0,0} ( X_{k,N}| )} \circ G(u) du  \\
= \int_0^{G^{-1 } ( \Vert  \E_{0,0} ( X_{k,\ell} ) \Vert_1)} Q_{ |\E_{0,0} ( X_{k,N}| )} (u) Q(u) du \, .
\end{multline*}
Using Inequality (4.6) in \cite{Rio17}, it follows that 
\[
\big \vert   \E \big ( \E_{0,0} ( X_{0,N} )  \E_{0,0} ( X_{k,\ell +N} )  \big )   \big  \vert 
\leq   \int_0^{ G^{-1 } ( \Vert  \E_{0,0} ( X_{k,\ell} ) \Vert_1)}  Q^2(u) du = \int_0^{ \Vert  \E_{0,0} ( X_{k,\ell} ) \Vert_1}  Q \circ G(u) du \, .
\]
Proceeding similarly, we have 
\[
\big \vert   \E \big ( \E_{0,0} ( X_{0,N} )  \E_{0,0} ( X_{k,\ell +N} )  \big )   \big  \vert  \leq  \int_0^{ \Vert  \E_{0,0} ( X_{0,N}   )\Vert_1}  Q \circ G(u) du \, 
\]
which goes to zero as $N \rightarrow \infty$ by condition \eqref{condTCLCF1}. Hence, by the dominated convergence theorem and condition \eqref{condTCLCF1} again, we derive that 
\[
\lim_{N \rightarrow \infty}  \sum_{k \geq 0}     \sum_{\ell \geq 0}    \big |  \E \big ( \E_{0,0} ( X_{0,N} )  \E_{0,0} ( X_{k,\ell +N} )  \big ) \big |  = 0 \, .
\]
With similar arguments, we get that 
\[
\lim_{N \rightarrow \infty}  \sum_{k \geq 0}     \sum_{\ell \geq 0}   \big |  \E \big ( \E_{0,0} ( X_{0,\ell +N} )  \E_{0,0} ( X_{k,N} )  \big ) \big | = 0 \, ,
\]
implying that 
\begin{equation} \label{neglireste2hN2}
\lim_{N \rightarrow \infty} \limsup_{n_1 \wedge n_2 \rightarrow \infty}\frac{1}{n_1n_2}   \Vert  S_{n_1,n_2} (h_{N,2})  \Vert_2^2 =0 \, ,
\end{equation}
 We proceed similarly to show that 
\begin{equation} \label{neglireste2}
\lim_{N \rightarrow \infty} \limsup_{n_1 \wedge n_2 \rightarrow \infty}\frac{1}{n_1n_2}  \big (  \Vert  S_{n_1,n_2} (h_{N,1})  \Vert_2^2 +   \Vert  S_{n_1,n_2} (h_{N,3})  \Vert_2^2 \big ) =0 \, ,
\end{equation}
where $h_{N,1} =\E_{0,0} ( U^N  f )$ and $h_{N,3} =\E_{0,0} (U^N V^N  f )$. The convergence \eqref{neglireste} (and then  \eqref{approxL2OMwithN})  follows from \eqref{neglireste2hN2} and \eqref{neglireste2}.  

\smallskip

Now notice that for $p$ and $q$ two positive  integers 
\begin{multline*}
 \Vert m_p - m_q \Vert_2 =  \frac{\Vert S_{n_1,n_2} (m_p)  - S_{n_1,n_2} (m_q)  \Vert_2}{\sqrt{ n_1 n_2}}  \\  \leq   \frac{\Vert S_{n_1,n_2} (f)  - S_{n_1,n_2} (m_p)  \Vert_2}{\sqrt{ n_1 n_2}}  +  \frac{\Vert S_{n_1,n_2} (f)  - S_{n_1,n_2} (m_q)  \Vert_2}{\sqrt{ n_1 n_2}}  \, .
\end{multline*}
Hence, the convergence  \eqref{approxL2OMwithN} implies that $(m_N)_{N \geq 1}$ is a Cauchy sequence in ${\mathbb L}^2$. We denote by $m$ its limit in ${\mathbb L}^2$. Then $(m\circ T_{{\bf{i}}})_{\bf{i} \in {\mathbb Z}^d}$ is a sequence of  orthomartingale differences  in ${\mathbb L}^2$. Next, we write
\[
 \frac{\Vert S_{n_1,n_2} (f)  - S_{n_1,n_2} (m)  \Vert_2}{\sqrt{ n_1 n_2}} \leq  \frac{\Vert S_{n_1,n_2} (f)  - S_{n_1,n_2} (m_N)  \Vert_2}{\sqrt{ n_1 n_2}} + \Vert m - m_N \Vert_2 \, ,
\]
which converges to zero by letting first $n_1 \wedge n_2$ tend to infinity and then $N$ and by considering \eqref{approxL2OMwithN}. This ends the proof of \eqref{approxL2OM}. 

\smallskip

Next, the CLT for $ \Big (  \frac{1}{\sqrt{n_1 n_2}} S_{n_1,n_2} (f)  \Big )  $  follows from the ${\mathbb L}^2$-orthomartingale approximation  \eqref{approxL2OM}  and the  CLT for stationary fields of  orthomartingales as proved in Voln\'y  \cite{V15} when at least one of the transformations generating the field is ergodic.  In addition, the fact that, under \eqref{condTCLCF1}, 
$ \displaystyle \lim _{n_1 \wedge n_2 \rightarrow \infty}\frac{1}{n_1n_2}  \Vert  S_{n_1,n_2} (f)  \Vert_2^2=\sigma^2 (f) $ follows by  standard arguments.  It remains to prove that 
$\Vert m \Vert_2^2 = \sigma^2(f)$.  This follows from the fact that, by  \eqref{approxL2OM}, 
\begin{multline*}
 \Big |  \Vert m \Vert_2  - \sigma^2 (f) \Big |  = \lim_{n_1 \wedge n_2 \rightarrow \infty} \Big |    \frac{\Vert S_{n_1,n_2} (m)   \Vert_2}{\sqrt{ n_1 n_2}}   -   \frac{\Vert S_{n_1,n_2} (f)   \Vert_2}{\sqrt{ n_1 n_2}} \Big |  \\
    \leq   \lim_{n_1 \wedge n_2 \rightarrow \infty} \frac{\Vert S_{n_1,n_2} (f)  - S_{n_1,n_2} (m)  \Vert_2}{\sqrt{ n_1 n_2}}  = 0  \, .
\end{multline*}
$\square$

\begin{remark}
Analyzing the proof of the theorem in case $d=2$, we infer that condition \eqref{condTCLCF1} can be replaced by the following set of conditions: 
\begin{equation} \label{condTCLCF12-bis}
\sum_{k, \ell \geq 0}    \sup_{u,v \geq 0} \big | \E \big (  \E_{0,0} ( X_{u,\ell} )  \E_{0,0} ( X_{k,v} )  \big )  \big | < \infty  \mbox{ and }\sum_{k, \ell \geq 0}    \sup_{u,v \geq 0} \big | \E \big (  \E_{0,0} ( X_{u,v} )  \E_{0,0} ( X_{k,\ell} )  \big )  \big | < \infty \, .
\end{equation}
Note that conditions \eqref{condTCLCF12-bis} are in the spirit of those given in Theorem 5.2 in \cite{HH80} in case $d=1$. 
\end{remark}


\subsection{The case of completely commuting transformations} \label{Section1.2}

Consider now a family $\{T_1, \ldots, T_d\}$ of measure preserving transformations on $(\Omega , {\mathcal A},  \mu )$. Denote by 
$U_1, \ldots, U_d$ the corresponding Koopman operators and by $K_1, \ldots, K_d$ the
associated adjoint operators. Those operators are characterized as follows
\[
U_i f = f \circ T_i
\]
 and 
 \begin{equation}\label{dual} {\mathbb E}( U_i f g ) =  {\mathbb E} (f K_i g)
\end{equation}
for every positive measurable functions $f,g$ and every $i \in \{1, \ldots, d \}$.  

\medskip

In particular, since $U_i(fg)=U_if U_ig$, we have ${\mathbb E}(  f g )={\mathbb E}( U_i f U_i g ) =  {\mathbb E} (f K_i U_ig)$ and we see that 
\begin{equation}\label{isometry}
K_iU_i=Id\, .
\end{equation}

\medskip

We shall assume that the family $\{T_1, \ldots, T_d\}$ (or the family $\{U_1, . . . , U_d\}$) is completely commuting in the sense of Gordin \cite{Go09} meaning that  it is commuting and  
\begin{equation} \label{propcom}
U_iK_j=K_jU_i \qquad \forall i\neq j \in \{1, \ldots , d\}\, .
\end{equation}

The natural filtrations associated with the transformations are defined as follows : for every $i \in \{1, \ldots, d \}$ and every $n \in {\mathbb N}$, ${\mathcal F}_n^{(i)} := T_i^{-n} ( \mathcal A ) $ and ${\mathcal F}_{{\bf{n}}} = \cap_{i=1}^d {\mathcal F}_{n_i}^{(i)}$.  It follows that for any $f \in {\mathbb L}^1 ( \Omega, {\mathcal A}, \mu )$, 
\[
\E ( f | {\mathcal F}_n^{(i)} ) =  U_i^n K_i^n f \, ,
\]
and, when $d=2$ (to simplify the exposition),
\begin{equation}\label{conditional-expectation}
\E(f |\F_{n_1,n_2})=\E\big( \E(f |\F^{(2)}_{n_2})|\F^{(1)}_{n_1}\big)= \E\big( \E(f |\F^{(1)}_{n_1})|\F^{(2)}_{n_2}\big)=
U_1^{n_1}U_2^{n_2} K_1^{n_1}K_2^{n_2} f \, .
\end{equation}
For $m$ in ${\mathbb L}^1( \mu)$ and ${\bf n} = (n_1, \dots, n_d )$, let $m_{{\bf n}} = U_1^{n_1}\cdots U_d^{n_d} m $.  We say that  $( m_{{\bf n}})_{{\bf n} \in {\mathbb Z}^d}$ is a sequence of reverse orthomartingale differences if $ \E ( m |  \F^{(i)}_{1} ) =0$ for any $i \in \{1, \dots, d \}$. 


%
%
%
%

\begin{theorem} \label{TCLCFreverse} Assume that $f$  is in ${\mathbb L}^2(\mu)$ and centered,  that 
 the family $\{T_1, \ldots, T_d\}$  is completely commuting and that at least one of the transformations is ergodic.  Let $S_{{\bf{n}}} (f) = \sum_{k_1=1}^{n_1} \ldots  \sum_{k_d=1}^{n_d}  U_1^{k_1}\cdots U_d^{k_d} f  $. Suppose in addition that 
\begin{equation} \label{condTCLCF1rev}
\sum_{\ell_1, \ldots \ell_d \geq 0}      \int_0^{ \min_{1 \leq i \leq d } \Vert  K^{\ell_i}_i f\Vert_1}  Q_{f} \circ  G_{f}  (u) du < \infty \, .
\end{equation}
Then,  there exists a sequence of  reverse orthomartingale differences $( m_{{\bf n}})_{{\bf n} \in {\mathbb Z}^d}$ in ${\mathbb L}^2 ( \mu)$ such that, as  $\min_{1 \leq i \leq d}n_i \ \rightarrow \infty$, 
\begin{equation*} \label{approxL2OMR}
 \frac{\Vert S_{{\bf{n}}} (f)  - S_{{\bf{n}}} (m)  \Vert_2}{\sqrt{|{\bf{n}}|}}  \rightarrow 0  \, .
\end{equation*}
Consequently, as  $\min_{1 \leq i \leq d}n_i \ \rightarrow \infty$, 
 \[  
\frac{1}{\sqrt{|{\bf{n}}|}} S_{{\bf{n}}} (f)   \rightarrow^{{\mathcal D}} {\mathcal N} ( 0, \Vert m \Vert_2^2)  \, .
\]
Moreover 
 \[  \Vert m \Vert_2^2 = \lim_{\min_{1 \leq i \leq d}n_i \ \rightarrow \infty}  \frac{\Vert S_{{\bf{n}}} (f)  \Vert_2^2}{|{\bf{n}}|} =  \sigma^2 (f)  \, ,
\]
where $\sigma^2 (f) =  \sum_{{\bf{k}} \in {\mathbb Z}^d} {\rm Cov}  (  X_{{\bf{k}}^+ } , X_{{\bf{k}}^- }) $ with 
${\bf{k}}^+ = (k_1^+, \ldots, k_d^+)$, ${\bf{k}}^- = (k_1^-, \ldots, k_d^-)$ and, for every ${\bf{\ell}} = ( \ell_1, \ldots, \ell_d) \in {\mathbb N}^d$, 
$X_{{\bf{\ell}}} = U_1^{\ell_1} \cdots  U_d^{\ell_d} f$. 
\end{theorem}
{\bf Proof of Theorem \ref{TCLCFreverse}.}  To simplify the exposition, we shall do the proof for $d=2$. For any $N $ a fixed positive integer, we first write the following decompositions: 
\[
f = \sum_{i=0}^{N-1} ( K_1^i f   -U_1 K_1^{i+1} f ) - (I- U_1) \sum_{i=0}^{N-1} K_1^{i+1} f + K_1^N f
\]
and
\[
f = \sum_{j=0}^{N-1} ( K_2^j f   -U_2 K_2^{j+1} f ) - (I- U_2) \sum_{j=0}^{N-1} K_1^{j+1} f + K_2^N f \, .
\]
This leads to the same decomposition as \eqref{co-martdecwithN}, that is : 
\begin{equation} \label{co-martdecwithN-rev}
  f = m_N - (I - U_1)g_{1,N} - (I - U_2)g_{2,N}  + (I-U_1)(I-U_2)g_{3,N} + h_N \, , 
\end{equation}
where 
\[
m_N = \sum_{i =0 }^{N-1} \sum_{j=0 }^{N-1}  (  K_1^i - U_1 K_1^{i+1}  ) (  K_2^j - U_2 K_2^{j+1}  ) f  \, , 
\]
\[
g_{1,N} = \sum_{i =0 }^{N-1} \sum_{j=0 }^{N-1} (  K_2^j - U_2 K_2^{j+1}  )  K_1^{i+1}  f \, , \, g_{2,N} = \sum_{i =0 }^{N-1} \sum_{j=0 }^{N-1}  (  K_1^i - U_1 K_1^{i+1}  )  K_2^{j+1} f   \, , 
\]
$g_{3,N} = \sum_{i =0 }^{N-1} \sum_{j=0 }^{N-1} K_1^{i+1} K_2^{j+1} f  $ and  
$h_N =   K_1^N f  + K_2^N f  - K_1^N K_2^N f  $.  We proceed now as in the proof of Theorem \ref{TCLCFd>2} using now the CLT for reverse orthomartingale differences (see Theorem 1 in \cite{CDV15}).  It remains to prove \eqref{approxL2OMwithN} which follows from \eqref{neglireste} with the same arguments as in the proof of Theorem  \ref{TCLCFd>2}.  To prove  \eqref{neglireste}, we  first prove that 
\begin{equation} \label{neglireste2hN2rev}
\lim_{N \rightarrow \infty} \limsup_{n_1 \wedge n_2 \rightarrow \infty}\frac{1}{n_1n_2}   \Big  \Vert  \sum_{i=1}^{n_1}  \sum_{j=1}^{n_2} U_1^i  U_2^j K_1^N f \Big \Vert_2^2 =0 \, .
\end{equation}
By the properties of the adjoint operators, 
\begin{multline*}
\Big \Vert  \sum_{i=1}^{n_1}  \sum_{j=1}^{n_2} U_1^i  U_2^j K_1^N f \Big \Vert_2^2 =
 4 \sum_{i=1}^{n_1} \sum_{k=i}^{n_1} \sum_{j=1}^{n_2} \sum_{\ell=j}^{n_2}    \E  ( K_1^{k-i +N} K_2^{\ell-j} (f  )  \cdot f )  \\
 + 4 \sum_{i=1}^{n_1} \sum_{k=i}^{n_1} \sum_{j=2}^{n_2} \sum_{\ell=1}^{j-1}    \E  ( K_1^{k-i +N} (f)  \cdot K_2^{j-\ell} (K_1^N f  ) ) \, .
\end{multline*}
Now, using \eqref{dual}, for $i \leq k$ and $\ell \leq j $, 
\begin{equation*}
\E  ( K_1^{k-i +N} (f)  \, K_2^{j-\ell} (K_1^N f  ) ) 
=  \E  (  f \, K_1^{k-i +N} f \circ T_1^N \circ T_2^{j-\ell}  )  \, .
\end{equation*}
By Proposition 1 in \cite{DD03} and  Inequality (4.6) in \cite{Rio17}, it follows that 
\[
| \E  ( K_1^{k-i +N} (f)  \cdot K_2^{j-\ell} (K_1^N f  ) ) |  \leq \int_0^{\Vert K_1^{k-i +N} (f) \Vert_1} Q_{|f|} \circ 
 G_{|f|} (u)  du \,.
 \]
 On another hand, proceeding similarly we get
 \[
| \E  ( K_1^{k-i +N} (f)  \cdot K_2^{j-\ell} (K_1^N f  ) ) |  \leq \int_0^{\Vert K_2^{j-\ell} (f) \Vert_1} Q_{|f|} \circ 
 G_{|f|} (u)  du \, ,
 \]
 giving
 \[
| \E  ( K_1^{k-i +N} (f)  \cdot K_2^{j-\ell} (K_1^N f  ) ) |  \leq \int_0^{ \min ( \Vert K_1^{k-i +N} (f) \Vert_1,  \Vert K_2^{j-\ell } (f) \Vert_1)} Q_{|f|} \circ 
 G_{|f|} (u)  du \, .
 \]
 A similar bound is valid for $
 | \E  ( K_1^{k-i +N} K_2^{\ell-j} (f  )  \cdot f ) |
$ when $\ell \geq j$.  By the dominated convergence theorem and condition \eqref{condTCLCF1rev},  the convergence \eqref{neglireste2hN2rev} follows. The rest of the proof of \eqref{neglireste} uses similar arguments. $\square$

\section{Examples} \label{Section2}

\setcounter{equation}{0}

\subsection{Linear random fields} \label{Section2RF}
For any $k , \ell$ in ${\mathbb Z}$, let 
\[
Y_{k, \ell}= \sum_{i,j \geq 0} a_{i,j} \varepsilon_{k-i,  \ell - j} \, ,
\]
where $|a_{i,j}| \leq  \kappa \rho^{i+j}$ for some $\rho \in (0,1)$ and $\kappa>0$,  and the random field $(\varepsilon_{i,j})_{i, j  \in {\mathbb Z}}$ 
is iid. Recall that, by Kolmogorov's three series theorem (see Theorem 5.3 in \cite{Kl14} for multiparameters r.v.'s),  we infer that 
this series is almost surely absolutely convergent as 
soon as 
\begin{equation}\label{moment}
{\mathbb E}( \log^2 (1 + |\varepsilon_{0,0}|))< \infty \, .
\end{equation}
Let $h$ be a bounded one Lipschitz-function and define 
\[
X_{k, \ell} = h( Y_{k, \ell}) - \E ( h( Y_{k, \ell}))  \text{ and } S_{n_1,n_2} = \sum_{k=1}^{n_1} \sum_{\ell=1}^{n_2} X_{k, \ell}  \, .
\]
\begin{corollary}\label{cor-LF-TCL}  If ${\mathbb E}( \log^4 (1 + |\varepsilon_{0,0}|))< \infty$,  the conclusion of Theorem \ref{TCLCFd>2} holds for 
$  \frac{S_{n_1,n_2}}{\sqrt{n_1n_2}}  $. 
\end{corollary}
{\bf Proof of Corollary \ref{cor-LF-TCL}.}
Let $(\varepsilon^*_{i,j})_{i, j  \in {\mathbb Z}}$ be an independent copy of $(\varepsilon_{i,j})_{i, j  \in {\mathbb Z}}$. For any non-negative integers $a$ and $b$, 
let $E_{a,b} = \{i \geq a, j \geq b\}$ and $S_{a,b} = {\mathbb N}^2 \backslash E_{a,b}$. Define 
\begin{equation} \label{notaA}
Y'_{k, \ell}  = \sum_{(i,j) \in E_{a,b} } a_{i,j} \varepsilon_{k-i,  \ell - j} +   \sum_{(i,j) \in S_{a,b} } a_{i,j} \varepsilon^*_{k-i,  \ell - j} 
\end{equation} 
and
\begin{equation} \label{notaB}
Y^*_{k, \ell} =  \sum_{i,j \geq 0} a_{i,j} \varepsilon^*_{k-i,  \ell - j} \, .
\end{equation} 
Clearly,
\begin{multline*}
 \big | \E_{k-a, \ell-b} ( X_{k,\ell} )  \big |  =   \big |  \E_{\varepsilon}  \big ( h(Y'_{k, \ell} ) - h(Y^*_{k, \ell} ) \big ) \big |  \\ \leq   \E_{\varepsilon} \Big ( 
 \min \Big ( 2 \Vert h \Vert_{\infty},  \sum_{(i,j) \in E_{a,b} } |a_{i,j} |  |  \varepsilon_{k-i,  \ell - j} -   \varepsilon^*_{k-i,  \ell - j} |   \Big ) \Big )  \, .
\end{multline*}
Hence, by  sub-additivity and stationarity, it follows that 
\begin{equation}  \label{notaC}
 \tau(a,b):= \big  \Vert  \E_{k-a, \ell-b} ( X_{k,\ell} ) \big \Vert_1 \leq  ( 2   \Vert h \Vert_{\infty} \vee 1 )   \sum_{(i,j) \in E_{a,b} }   \E \big ( \min  (  |a_{i,j} |  |  \varepsilon_{0,0} -   \varepsilon^*_{0, 0} | , 1  )  \big ) \, .
\end{equation} 
Now, letting $c=\rho^{-1/2}$, we have 
\begin{equation*}
 \tau(a,b) \leq 
\kappa \sum_{i \geq a} \sum_{j \geq b} \rho^{i+j} {\mathbb E}\left(  |\varepsilon_0-\varepsilon^*_0| {\bf{1}}_{|\varepsilon_0-\varepsilon^*_0|\leq c^{i+j}} \right ) 
+\sum_{i \geq a} \sum_{j \geq b}
{\mathbb P}\left ( |\varepsilon_0-\varepsilon^*_0| > c^{i+j}\right )
\, .
\end{equation*}
Hence 
\[
 \tau(a,b) \leq  
 \frac{\kappa}{(1 - \sqrt{\rho})^2} \rho^{(a+b)/2} + {\mathbb E}\left 
( \left ( \frac{\log |\varepsilon_0-\varepsilon^*_0|}{\log(c)}
-a\right)_+ \left ( \frac{\log |\varepsilon_0-\varepsilon^*_0|}{\log(c)}
-b\right)_+ \right) \, , 
\]
and finally
\[
 \tau(a,b) \leq  K \rho^{(a+b)/2}+ K{\mathbb E}\left 
(  ( \log |\varepsilon_0-\varepsilon^*_0| )^2
{ \bf 1}_{\log |\varepsilon_0-\varepsilon^*_0|> (a \vee b) \log(c)}\right) \, ,
\]
for some $K>0$.
Note that $ \tau(a,b)  \rightarrow 0$ as $a \vee b  \rightarrow \infty$
as soon as \eqref{moment} is satisfied.

Recall that $\gamma_1(k) = \Vert \E_{0,0} ( X_{k,0} ) \Vert_1$ and $\gamma_2(k) = \Vert \E_{0,0} ( X_{0,k} ) \Vert_1$.  Setting 
\begin{equation} \label{deflambdak}
\lambda(k) :=  K \rho^{k/2}+ K{\mathbb E}\left 
(  ( \log |\varepsilon_0-\varepsilon^*_0| )^2
{ \bf 1}_{\log |\varepsilon_0-\varepsilon^*_0|> k \log(c)}\right)  \, , 
\end{equation}
we have for $i=1,2$, 
\[
\gamma_i(k) \leq \lambda(k) \, .
\]
Consequently, condition \eqref{condTCLCF1} holds as soon as $\sum_{k \geq 1} k \lambda(k) < \infty$, which is satisfied as son as ${\mathbb E}\left (\log^{4} (1 + |\varepsilon_{0,0}|)\right) 
< \infty$. $\square$

\subsection{Functions of expanding endomorphisms of the $m$-dimensional torus} \label{SectionEE}

In this section, we shall apply Theorem \ref{TCLCFreverse} to the case where the transformations  $(T_1, \ldots , T_d)$ are given by expanding endomorphisms of the $m$-dimensional torus $\TT_m:=\R^m/\Z^m$.  We use the setup of Section 5 in \cite{CDV15}.  

\medskip

An expanding endomorphism of $\TT_m$ is of the form 
$T_A\, :\, x\mapsto Ax$, where $A$ is an $m\times m$ matrix with integral entries and all eigenvalues of modulus strictly greater than 1. We shall refer to \cite{CDV15} for all the needed properties 
of $T_A$, some of those are recalled below.

\medskip

Recall that $T_A$ preserves the Haar-Lebesgue measure $\lambda_m$ on $\TT_m$. We then denote by $U_A$ the Koopman operator associated with $T_A$. 

\medskip

Recall that the ${\mathbb L}^1$-modulus of continuity $\omega_{1,f}$ of $f
\in {\mathbb L}^1(\TT_m,\lambda_m)$ is given by 
$$
\omega_{1,f}(\delta) =\sup_{h\in  \TT_m\, :\, |h|\le \delta } 
\|f-f(\cdot + h)\|_1\, .
$$ 
with $|h|=\inf_{\ell \in \Z^m}|{\tilde h}+\ell|_2$ where $ {\tilde h}$ is the representative of $h$ in   $ [0,1)^m$ and $|\cdot |_2$ stands for the euclidean norm on $\R^m$.

Recall the following lemma which  gives a simple condition under which two expanding endomorphisms are completely commuting.
\begin{lemma}[Cuny-Dedecker-Voln\'y \cite{CDV15}] Let $A$ and $B$ be two expanding $m \times m$ ($m \geq 1$) matrices with integer entries.
Assume that $A$ and $B$ commute and that they have coprime determinants ($\gcd(det A , det B)$ = 1). Then, $T_A$ and $T_B$ are completely commuting.
\end{lemma}

\medskip

Let us now state the main result of this section.

\begin{corollary}\label{theo-transfo}
Let $T_{A_1}, \ldots, T_{A_d}$ be  expanding endomorphisms of $\TT_m$ that are completely commuting. Let $f\in  L^\infty(\TT_m,\lambda_m)$ and centered. Assume that 
\begin{equation} \label{condmodcont}
\int_0^1 \frac{|\log (t) |^{d-1}}{t} \omega_{1,f} (t) dt  < \infty 
\end{equation}
Then the conclusion of Theorem \ref{TCLCFreverse} holds. 
\end{corollary}
Theorem 11 in \cite{CDV15} asserts that the conclusion of Theorem \ref{TCLCFreverse} holds if $f$ is in $L^2(\TT_m,\lambda_m)$, centered  and such that 
\begin{equation} \label{condmodcontL2}
\int_0^1 \frac{|\log (t) |^{(d-2)/2}}{t} \omega_{2,f} (t) dt  < \infty  \, , 
\end{equation}
where $
\omega_{2,f}(\delta) =\sup_{h\in  \TT_m\, :\, |h|\le \delta } 
\|f-f(\cdot + h)\|_2
$. Now when $f$ is in ${\mathbb L}^{\infty}(\TT_m,\lambda_m)$, one can control $\omega_{2,f}$ by $\omega_{1,f}$ via the inequality
$
\omega_{2,f}(\delta) \leq \sqrt{ 2 \Vert f \Vert_{\infty} \omega_{1,f}(\delta)  } $. Hence, in the bounded case, the condition
\begin{equation} \label{condmodcontL1bis}
\int_0^1 \frac{|\log (t) |^{(d-2)/2}}{t}  \sqrt{\omega_{1,f} (t)} dt  < \infty  \, , 
\end{equation}
implies condition \eqref{condmodcontL2}.  Note that condition \eqref{condmodcontL1bis} implies condition \eqref{condmodcont}.  Therefore, in the bounded case, if we are interested by criteria involving the ${\mathbb L}^1$-modulus of continuity it is preferable to use Corollary \ref{theo-transfo} rather than condition \eqref{condmodcontL1bis}. 

\begin{example} In case $m=1$, as an application we can consider a family of transformations $(T_{p_i})_{1 \leq i \leq d}$, defined, for any $x \in [0,1]$, by
\[
T_{p_i} (x) = p_i x  -  [p_i x]  \, , 
\]
where $(p_i)_{1 \leq i \leq d}$ is a family of integers which are pairwise coprime. For such maps,  we can consider observables of the type $f = g ( \varphi) $ where $g$ is Lipshitz and bounded from ${\mathbb R}$ to ${\mathbb R}$ and  $\varphi$ is defined as page 118 in \cite{BLT}. More precisely, let $\omega$ be a concave modulus of continuity satisfying \eqref{condmodcont}.  We then take 
\[
\varphi (x) = \sum_{k \geq 1} \mu^* (k) \cos (kx) 
\]
where $\mu^* (t)$ is the greatest convex minorant of $\mu (t) = \omega(t^{-2})$.  Then $\omega_{1,f} (\cdot)$ satisfies 
\eqref{condmodcont}.  Note that the modulus of continuity of $\varphi$ can be evaluated in ${\mathbb L}^1$ but not in ${\mathbb L}^2$ since its Fourier coefficients are not in ${\mathbb L}^2$ in the interesting cases. Hence it seems difficult to give for $\omega_{2,f} ( \delta )$ a better control than $
\omega_{2,f}(\delta) \leq \sqrt{ 2 \Vert f \Vert_{\infty} \omega_{1,f}(\delta)  } $.
\end{example}

\medskip

Before doing the proof of Corollary \ref{theo-transfo}, we shall provide some notations and facts 
about expanding endomorphisms of $\TT_m$. 

\medskip

If $\Gamma \subset \Z^m$ is a system of representatives of $\Z^m/A\Z^m$, then the adjoint operator $K_A$ of $U_A$ is given, for every positive measurable $f$, by 

\begin{equation}\label{perron}
K_Af(x)=\frac1{{\rm det} \, A} \sum_{\gamma\in \Gamma} f(A^{-1}x+A^{-1}\gamma) \quad \forall x\in \TT_m\, .
\end{equation}
\medskip

Given a matrix $M\in M_m(\R)$ we denote by $|M|_2$ its operator 
norm associated with the euclidean norm. 

\medskip

The following lemma is just the ${\mathbb L}^1$-version of Proposition 13 of \cite{CDV15}. The proof is identical. 

\begin{lemma}\label{lemma-transfo}
For every $f\in {\mathbb L}^1(\TT_m,\lambda_m)$ with $\lambda_m(f)=0$, we have 
\[
\| K_Af\|_1\le \omega_{1,f}\big(\Delta(\, A^{-1}([0,1)^m)\,)\big) \, , \]
 where 
$\Delta\big(A^{-1}([0,1)^m)\big)=\sup_{x\in [0,1)^m}|A^{-1}x|$.
\end{lemma} 

We shall also need the following.

\begin{lemma}\label{lemma-transfo-2}
Let $A$  be an expanding $m \times m$ ($m \geq 1$) matrix with integer entries.   For every $f\in {\mathbb L}^1(\TT_m,\lambda_m)$, 
we have 
\begin{equation}\label{modulus-K}
\omega_{1, K_A f}(\delta)\le \omega_{1,f} (|A^{-1}|_2\delta)
\quad \forall \delta>0\, .
\end{equation}
\end{lemma}
\noindent {\bf Proof of Lemma \ref{lemma-transfo-2}.} 
Let $\delta >0$. Let $h \in \TT_m$  be such that $| h|\le \delta$, and let ${\tilde h}$ be its representative in $[0,1)^m$. There exists $\ell\in \Z^m$ such that $ | {\tilde h}+\ell |_2 = |h|\le \delta$. 

\medskip

Set, for every $x\in \TT_m$, $\varphi(x):=  f(x) -f(x+A^{-1}(h+\ell))\, .$

\medskip

Using \eqref{perron} and the fact that $K_A f(\cdot +h)=K_A f(\cdot +h+\ell)$, we see that 
$$
\|K_A f - K_Af(\cdot +h)\|_1=\|K_A f - K_Af(\cdot +h+\ell)\|_1= \|K_A \varphi\|_1\le \|\varphi \|_1\, .
$$
Next note that 
$$
|A^{-1}(h+\ell) | \le |A^{-1}|_2 |{\tilde h}+ \ell|_2 = |A^{-1}|_2 |h|\le  |A^{-1}|_2  \delta\, ,
$$
and the result follows.  \hfill $\square$

\medskip

Finally, we recall without proof the following standard lemma. 

\begin{lemma}\label{lemma-transfo-3}
For every $f,g\in {\mathbb L}^1(\TT_m,\lambda_m)\cap L^\infty(\TT_m,\lambda_m)$ we have 
\begin{equation}
\omega_{1,fg}(\delta)\le \|f\|_\infty \omega_{1,g}(\delta) +\|g\|_\infty \omega_{1,f}(\delta)\quad \forall \delta >0\, .
\end{equation}
\end{lemma}

We are in position now to give the proof of Corollary \ref{theo-transfo}.

\medskip

\noindent {\bf Proof of Corollary \ref{theo-transfo}.}  To simplify the exposition, let us consider the case $d=2$.  Since $f$ is bounded, the function $Q_f$ is bounded  and we are back to prove that 
\begin{equation} \label{series-cond-0}
\sum_{k, \ell \geq 0}        \Vert  K_{A_1}^k(f) \Vert_1 \wedge  \Vert  K_{A_2}^\ell(f) )  \Vert_1 < \infty \, .
\end{equation}

Now, notice  that, for every $i\in \{1,2\}$ and every integer $n\in \N$, $K_{A_i}^n=K_{A_i^n}$. Let $\eta$ be the maximum of the modulus of the eigenvalues of $A_1$ and $A_2$. By assumptions, 
$\eta<1$. Hence there exists $C>0$ such that for every $n\in \N$, we have, for every $i\in \{1,2\}$
\begin{equation}\label{contraction}
\Delta\big(A_i^{-n}([0,1)^m)\big)\le C \eta ^n\, .
\end{equation}
Hence, in view of Lemma \ref{lemma-transfo-3}, we have 
\begin{equation*}
\Vert  K_{A_1}^k(f) \Vert_1 \wedge  \Vert  K_{A_2}^\ell(f) )  \Vert_1\le \omega_{1,f}(C\eta^{\max (k,\ell)})\, .
\end{equation*}

In particular, \eqref{series-cond-0} holds as soon as 
\begin{equation}\label{series-cond}
\sum_{n\in \N} n \omega_{1,f}(C\eta^n)<\infty\, .
\end{equation}
Comparing series with integrals and using that $\omega_{1,f}$ is nondecreasing, we infer that \eqref{series-cond} is equivalent to 
$$\int_0^\infty x\,  \omega_{1,f}(\delta^x)dx <\infty
$$
which, by a change of variable, is equivalent to  condition \eqref{condmodcont} with $d=2$. This ends the proof of the corollary. $\square$

\section{The weak invariance principle} \label{section3}
\setcounter{equation}{0}


\subsection{Weak dependence conditions for commuting filtrations} \label{Section2.1}

We are in the setup of Section \ref{Section1.1}.  Let us introduce the following weak dependence coefficients  that are well adapted to prove the functional part of the central limit theorem. For any positive real $M$, let $\varphi_M(x) = ( x \wedge M) \vee (-M)$. We assume that there exists a non-increasing  sequence $(\theta (k) )_{k \geq 0} \subset {\mathbb R}^+$ such that 
\[
\sup_{M>0}  \Vert \E_{{\bf{0}}} ( \varphi_M(X_{{\bf{i}}} )  ) -  \E( \varphi_M(X_{{\bf{i}}} )  )  \Vert_1     \leq  \min_{1 \leq \ell \leq d }  \theta ( i_{\ell} )  
\]
and
\[
\sup_{M>0} M^{-1} \Vert \E_{{\bf{0}}} ( \varphi_M(X_{{\bf{i}}} ) \varphi_M (X_{{\bf{j}}} ) ) -  \E( \varphi_M(X_{{\bf{i}}} ) \varphi_M (X_{{\bf{j}}} ) )  \Vert_1  \leq   \min_{1 \leq \ell \leq d }  \theta ( i_{\ell} \wedge j_{\ell})  \, .
\]

When the random field is bounded by $M$, we can also use the following weak dependent coefficients. Let $(\gamma (k) )_{k \geq 0} \subset {\mathbb R}^+$ be a 
non-increasing sequence sequence  such that 
\[
\Vert \E_{{\bf{0}}} (X_{{\bf{i}}} )   \Vert_1     \leq  \min_{1 \leq \ell \leq d }  \gamma ( i_{\ell} )  
\]
and
\[
 \Vert \E_{{\bf{0}}} ( X_{{\bf{i}}} X_{{\bf{j}}}  ) -  \E( X_{{\bf{i}}} X_{{\bf{j}}}  )  \Vert_1  \leq   M  \min_{1 \leq \ell \leq d }  \gamma ( i_{\ell} \wedge j_{\ell})  \, .
\]
For every $ {\bf  n}=(n_1, \ldots,n_d) \in {\mathbb N}^d$ and $ {\bf  t}=(t_1, \ldots,t_d) \in [0,1]^d $, set
\[
S_{ {\bf  n},  {\bf  t}} (f) := \sum_{ {\bf  0} \preceq {\bf k} \prec  [ {\bf  n t }]}  \prod_{i=1}^d (k_i \wedge (n_it_i -1) -k_i +1 )U_{ {\bf  k}} f  \, , 
\]
where $[ {\bf  n t }]= ( [n_1t_1], \ldots, [n_dt_d]) $. Define also
\begin{equation} \label{defdeW}
W_{ {\bf  n},  {\bf  t}} (f) :=  \frac{S_{ {\bf  n},  {\bf  t}} (f) }{  \big ( \prod_{i=1}^d n_i \big )^{1/2}} \, .
\end{equation}

\begin{theorem} \label{FCLTthm}
Assume one of the following conditions : 
\begin{enumerate}
\item[1.] $\Vert f \Vert_{\infty} \leq M$ for some $M >0$, and $ \sum_{k \geq 1} k^{d-1} \gamma(k) < \infty$.
\item[2.] $\sum_{k \geq 0} (k+1)^{d-1}  \int_0^{\theta(k)} Q_f \circ G_f (u) du  < \infty$.
\end{enumerate}
Then, when $\min_{1 \leq i \leq d } n_i \rightarrow \infty$, the sequence of processes $( \{W_{ {\bf  n},  {\bf  t}} (f) ,  {\bf  t} \in [0,1]^d \} )_{ {\bf  n} \in {\mathbb N}^d}$ converges in law in $( C([0,1]^d ) , \Vert \cdot \Vert_{\infty} )$ to 
$ \{ \sigma (f) W_{  {\bf  t}} ,  {\bf  t} \in [0,1]^d  \}$ where $ \{ W_{  {\bf  t}} ,  {\bf  t} \in [0,1]^d  \}$ is the standard $d$-dimensional brownian sheet and 
$\sigma^2 (f) =  \sum_{ {\bf  k} \in {\mathbb Z}^d} {\rm Cov} ( X_ {\bf  0},  X_ {\bf  k} )$. 
\end{theorem}
Condition of Item 2. can be rewritten as 
\begin{equation} \label{condtightness} 
\int_0^1 Q_f \circ G_f(u)   ( \theta^{-1} (u) )^d du < \infty \, , 
\end{equation}
which in turn is equivalent to 
\[
 \int_0^1 R_f(u) Q_f(u)  du < \infty  \text{ where } R_f(u)  = Q_f(u)  ( \theta^{-1} \circ G_f^{-1} (u) )^{d}\, . 
\]

\subsection{The case of completely commuting transformations}

We are in the setup of Section \ref{Section1.2} assuming in addition that  $f$ is  bounded to simplify the exposition.  We shall consider the following dependence condition: there exists some decreasing function $\gamma$ such that for every $ {\bf  k}, {\bf \ell} , {\bf i} \in \N^d$,
\begin{equation}\label{plusieurs-points-bis}
\|\E(f|\F_{ {\bf i}})\|_1 \le \min_{1 \leq j \leq d } \gamma(i_j)\, ,
\end{equation}
\begin{equation}\label{plusieurs-points}
\|\E(U_1^{k_1}\ldots U_d^{k_d}f U_1^{\ell_1}\ldots U_d^{\ell_d}  f|\F_{ {\bf k} \vee  {\bf \ell}  +  {\bf i}})-
\E(U_1^{k_1}\ldots U_d^{k_d}f U_1^{\ell_1}\ldots U_d^{\ell_d}  f)\|_1\le  \Vert f \Vert_{\infty} \min_{1 \leq j \leq d } \gamma(i_j)\, ,
\end{equation}
and 
\begin{equation}
\label{condplusieurs-points} 
\sum_{k\ge 1}k^{d-1}\gamma(k)<\infty\, .        
\end{equation}
\begin{theorem} \label{FCLTthmreverse} Assume that $f$  is bounded and centered, that 
 the family $\{T_1, \ldots, T_d\}$  is completely commuting and that at least one of the transformations is ergodic.  Let $W_{ {\bf  n},  {\bf  t}} (f)$ be defined in  \eqref{defdeW} with $U_{{\bf i}}  = U_1^{i_1}\ldots U_d^{i_d}$.  Assume in addition that condition \eqref{condplusieurs-points} holds. Then, when $\min_{1 \leq i \leq d } n_i \rightarrow \infty$, the sequence of processes $( \{W_{{\bf n}, {\bf  t}} (f) , {\bf  t} \in [0,1]^d  \} )_{ {\bf  n} \in {\mathbb N}^d}$ converges in law in  $( C([0,1]^d ) , \Vert \cdot \Vert_{\infty} )$ to 
$ \{ \sigma (f) W_{ {\bf  t}} , {\bf  t} \in [0,1]^d  \}$ where $ \{ W_{ {\bf  t}} , {\bf  t} \in [0,1]^d  \}$ is the standard $d$-dimensional brownian sheet and 
$\sigma^2 (f) $ is defined in Theorem \ref{TCLCFreverse}. 
\end{theorem}

\subsection{Examples (continued)} \label{examplescontinued}

We consider the examples of Section \ref{Section2}.

\subsubsection{Linear random fields}

Let us consider the linear random fields as defined in Section \ref{Section2RF}.  Under the same condition as in Corollary \ref{cor-LF-TCL}, the functional form of the CLT holds. 

\begin{corollary}\label{cor-LF-TCLfonctionnel}  If ${\mathbb E}( \log^4 (1 + |\varepsilon_{0,0}|))< \infty$,  the functional form of Corollary \ref{cor-LF-TCL}  holds.  
\end{corollary}
{\bf Proof of Corollary \ref{cor-LF-TCLfonctionnel}.}  We verify the conditions of Theorem \ref{FCLTthm}. 
Clearly, with the notations \eqref{notaA} and \eqref{notaB}, 
\begin{multline*}
 \big |  \E_{0,0} ( h(Y_{i,j} ) h( Y_{k,\ell}  )) -  \E ( h(Y_{i,j} ) h( Y_{k,\ell}  )) \big |  =   \big |  \E_{\varepsilon}  \big ( h(Y'_{i,j} ) h(Y'_{k, \ell} ) -h(Y^*_{i,j} )  h(Y^*_{k, \ell} ) \big ) \big |  \\ \leq  \Vert h \Vert_{\infty}  \big |  \E_{\varepsilon}  \big ( h(Y'_{i,j} )  -h(Y^*_{i,j} )  \big ) \big | +    \Vert h \Vert_{\infty}  \big |  \E_{\varepsilon}  \big ( h(Y'_{k,\ell} )  -h(Y^*_{k,\ell} )  \big ) \big |  \, .
\end{multline*}
Therefore, with the notation \eqref{notaC}, 
\[
\Vert   \E_{0,0} ( h(Y_{i,j} ) h( Y_{k,\ell}  )) -  \E ( h(Y_{i,j} ) h( Y_{k,\ell}  )) \Vert_1 \leq    \Vert h \Vert_{\infty}  \big ( \tau(i,j) + \tau(k, \ell) \big )  \, .
\]
Hence, with $\lambda(k)$ defined by \eqref{deflambdak}, we get 
\[
\Vert   \E_{0,0} ( h(Y_{i,j} ) h( Y_{k,\ell}  )) -  \E ( h(Y_{i,j} ) h( Y_{k,\ell}  )) \Vert_1 
 \leq K \min ( \lambda ( k \wedge i ) , \lambda ( \ell \wedge j )  ) \, .
\]
In addition
\[
\Vert   \E_{0,0} ( h(Y_{i,j} ) ) -  \E ( h(Y_{i,j} )  \Vert_1 
 \leq K \min ( \lambda ( i ) , \lambda (  j )  ) \, .
\]
So, the condition of Theorem \ref{FCLTthm} is satisfied as soon as $ \E (  \log^4 ( 1+  |\varepsilon_{0,0}| )  ) < \infty$.  $\square$

\subsection{Functions of expanding endomorphisms of the $m$-dimensional torus}

Let us consider the commuting endomorphisms as in  Section \ref{SectionEE}.  Under the same conditions as in Corollary \ref{theo-transfo}, the functional form of the CLT holds.

\begin{corollary}\label{theo-transfo-FCLT}
Let $T_{A_1}, \ldots, T_{A_d}$ be  expanding endomorphisms of $\TT_m$ that are completely commuting. Let $f\in  L^\infty(\TT_m,\lambda_m)$ and centered. Assume that  condition \eqref{condmodcont} holds.  Then the conclusion of Theorem \ref{FCLTthmreverse}  holds.  
\end{corollary}
{\bf Proof of Corollary \ref{theo-transfo-FCLT}.} Once again to simplify the exposition, we consider the case $d=2$. The result follows by an application of Theorem \ref{FCLTthmreverse}  provided one can verify \eqref{plusieurs-points} and \eqref{condplusieurs-points}.

Let us prove  \eqref{plusieurs-points}. 
 Let $k,\ell,r,s\in \N_0=\{0,1,\ldots\}$. Set $u:= K_{A_1}^kK_{A_2}^\ell f K_{A_1}^rK_{A_2}^sf$. By Lemma \ref{lemma-transfo-2} and Lemma \ref{lemma-transfo-3}, 
$$
\omega_{1,u}(\delta)\le 2\|f\|_\infty \omega_{1,f}(D\delta),
$$
with $D:=\sup_{m,n\in \N_0} |A_1^{-m}|_2|A_{2}^{-n}|_2<\infty$. Using Lemma \ref{conditional-lemma}, Lemma 
\ref{lemma-transfo} and the fact that $K_{A_1}$ and $K_{A_2}$  are contractions of ${\mathbb L}^1(\TT_m,\lambda_m)$, we infer that for every $k,\ell,r,s, j,q\in \N$,
\begin{gather*}
\|\E(U_1^kU_2^rfU_1^\ell U_2^s f|\F_{k\vee \ell+j,r\vee s+q})-\E(
U_1^kU_2^rfU_1^\ell U_2^s f)\|_1 \qquad \qquad\\ \qquad 
\qquad
\le 2\|f\|_\infty \min \big(\omega_{1,f} \big( D\Delta (A_1^{-j}(|0,1)^m))\big)
, \omega_{1,f} \big( D\Delta (A_2^{-q}(|0,1)^m))\big)\big)\, ,
\end{gather*}
Hence, using \eqref{contraction}, the first part of  \eqref{plusieurs-points} holds with $\gamma$ defined by
$$
\gamma(k):= 2\|f\|_\infty \omega_{1,f}(CD\eta ^k),\quad k\in \N_0\, .
$$
We can obtain a similar bound for \eqref{plusieurs-points-bis}. In addition  we have already seen that the condition $\sum_{k\ge 0} 
k\omega_{1,f}(CD\eta ^k)<\infty$ holds under our assumption, which proves \eqref{condplusieurs-points}. \hfill $\square$


\section{Proofs of the results of Section \ref{section3}}

In this section and the appendix, for any two positive sequences of random variables $a_{n}$
and $b_{n}$ we shall often use the notation  $a_{n}\ll b_{n}$  to mean that there is a 
positive finite constant $c$  such that $a_{n}\leq cb_{n}$ for all $n$. 

\setcounter{equation}{0}

\subsection{Proof of Theorem \ref{FCLTthm}}

The result will follow if one can prove the convergence of the finite dimensional laws and the tightness of the process.

\medskip

\noindent {\it 1. Convergence of the finite dimensional laws.} To simplify the exposition, let us consider the case $d=2$ (the general case being identical but uses more notations). Set $T_{1,0}=T$ and $T_{0,1}=S$ and assume that the  transformations $T$ and $S$ are commuting and that $T$ (for instance) is ergodic. Let $U$ and $V$ the operators defined by 
$Uf = f \circ T$ and $Vf = f \circ S$.
To prove the convergence of the finite dimensional laws, it suffices to show that for any $K \geq 1$, any $0=t_0 < t_1 <\ldots < t_K =1 $ and $0=s_0 < s_1 <\ldots < s_K =1 $, and any reals $(a_{k, \ell})_{1 \leq k, \ell \leq K}$, 
\[
V_{n_1,n_2} (f)  = \frac{1}{\sqrt{n_1n_2}} \sum_{k, \ell = 1}^K a_{k , \ell} \sum_{i=[n_1 t_{k-1}]  +1}^{[n_1 t_k] } \sum_{j=[n_2 s_{\ell-1}]  +1}^{[n_2 s_\ell] }  U^i  V^j f 
\]
converges in law, as $n_1,n_2 \rightarrow \infty$ to
\[
 \sigma( f) \sum_{k, \ell = 1}^K a_{k , \ell} \big ( W_{t_k, s_\ell} + W_{t_{k-1}, s_{\ell-1}} - W_{t_k, s_{\ell-1}} - W_{t_{k-1}, s_{\ell}}   \big ) \, .
\]
(See for instance the explanations given Section 3.2 in \cite{CDV15}.)  To prove the above convergence, we proceed as in the proof of  Theorem  \ref{TCLCFd>2}:  for a positive integer $N$, we 
consider the decomposition \eqref{co-martdecwithN}, namely
\[
  f = m_N + (I - U)g_{1,N} + (I - V)g_{2,N}  + (I-U)(I-V)g_{3,N} + h_N \, . 
\]
With the same arguments as those developed in the proof of Theorem  \ref{TCLCFd>2}, we get  
\begin{equation*} \label{AMFDL}
 \frac{\Vert  V_{n_1,n_2} (f) - V_{n_1,n_2} (m)   \Vert_2}{\sqrt{ n_1 n_2}}  \rightarrow 0  \, , \, \text{as  }n_1 \wedge n_2  \rightarrow \infty \, , 
\end{equation*}
where $m$ is the limit in ${\mathbb L}^2$ of $m_N$. Using the convergence (7) in \cite{CDV15}, we obtain that 
\[
 V_{n_1,n_2} (m)  \rightarrow^{{\mathcal D}} G \, ,  \, \text{ as $n_1 \wedge n_2 \rightarrow \infty$}, 
\]
where $G \sim {\mathcal N} (0,  \Gamma)$ with 
\[
\Gamma:= \Vert m \Vert^2_2  \sum_{k, \ell = 1}^K a^2_{k , \ell} (t_k - t_{k-1} ) ( s_{\ell} - s_{\ell-1} ) \, .
\]
This ends the proof of  the convergence of the finite dimensional laws.

\medskip

\noindent {\it 2. Tightness.} The next result proves the tightness of the partial sums process. 
\begin{proposition} \label{Proptight}
Under the conditions of Theorem \ref{FCLTthm}, 
 \begin{equation} \label{tightaim}
\lim_{ \lambda \to \infty} \limsup_{\min_{1 \leq \ell \leq d } n_{\ell} \to \infty}  \lambda^2 
  \p \Big ( \max_{{\bf{i}} \preceq {\bf{n}}} |S_{{\bf{i}}} (f) | \geq \lambda \sqrt{|{\bf{n}}| } \Big ) = 0 \, . 
\end{equation} 
\end{proposition}
{\bf Proof of Proposition \ref{Proptight}.}  We shall prove the proposition under condition \eqref{condtightness}. In case when $\Vert f \Vert_{\infty} \leq M$, the proof will use the same arguments with some simplifications since we do not need to truncate the random variables. To soothe the notations we write $Q, G, R$ instead of $Q_f, G_f, R_f$. 

 To simplify the exposition, we do the proof in case $d=2$ and we set $T_{1,0}=T$ and $T_{0,1}=S$. We  assume that $T$ and $S$ are commuting and we denote by $U$ and $V$ the operators defined by 
$Uf = f \circ T$ and $Vf = f \circ S$.  Let $v \in [0,1]$ and $M=Q(v)$. We set 
\[
X'_{k,\ell} = \varphi_M (X_{k,\ell} )- \E (\varphi_M (X_{k,\ell} )) \mbox{ and }  X''_{k,\ell} = X_{k,\ell} - X'_{k,\ell} \, .
\]
Let $S'_{i,j} (f) = \sum_{k=1}^i \sum_{\ell=1}^j X'_{k,\ell}$ and $S''_{i,j} (f) = \sum_{k=1}^i \sum_{\ell=1}^j X''_{k,\ell}$.  Observe also that for $M = Q(v)$ with $v \in [0,1]$, 
\[
Q_{|\varphi_M (X_{0,0} )|} (u) \leq  Q( u \vee v) \text{ and } Q_{|X_{0,0} -\varphi_M (X_{0,0} )|} (u) \leq   Q( u )  {\bf 1}_{u \leq v} \, .
\]
We shall choose $v$ as follows. Let 
\[R = Q (\theta^{-1} \circ G^{-1})^2 ,  \, x = \sqrt{n_1n_2}  \text{ and  } v = R^{-1} ( x ) \, . \] 
Write first
\begin{multline*} 
  \p \Big (\max_{i \leq n_1 \atop j \leq n_2} |S_{i,j}(f)|  \geq 6 \lambda \sqrt{n_1 n_2 } \Big ) \\
   \leq  \p \Big (\max_{i \leq n_1 \atop j \leq n_2} |S'_{i,j}(f)|  \geq 5 \lambda \sqrt{n_1 n_2 } \Big )  + \p \Big (\max_{i \leq n_1 \atop j \leq n_2} |S''_{i,j}(f)|  \geq  \lambda \sqrt{n_1 n_2 } \Big )  \, .
\end{multline*}
Using Markov's inequality and stationarity, we get  
\begin{multline*} 
  \p \Big (\max_{i \leq n_1 \atop j \leq n_2} |S''_{i,j}(f)|  \geq  \lambda \sqrt{n_1 n_2 } \Big )  \leq   \frac{2  \sqrt{n_1 n_2 } }{ \lambda }  
   \E (|X_{0,0} -\varphi_M (X_{0,0} )|)  \\
   \leq    \frac{2  \sqrt{n_1 n_2 } }{ \lambda }  
  \int_0^1 Q(u)  {\bf 1}_{u \leq  R^{-1} ( x ) }  du \, .
\end{multline*}
Since $R$ is right-continuous and non-increasing,
\begin{equation} \label{Rrightcontinuous}
u < R^{-1} ( x )  \iff R(u) > x \, .
\end{equation}
Hence, since $x = \sqrt{n_1 n_2}$,
 \begin{equation} \label{tightaimtruemoinsboundedpart}
  \p \Big (\max_{i \leq n_1 \atop j \leq n_2} |S''_{i,j}(f)|  \geq  \lambda \sqrt{n_1 n_2 } \Big )  \leq   \frac{2   }{ \lambda }  
  \int_0^1 R(u) Q(u)  {\bf 1}_{R(u) \geq  x }  du \, ,
\end{equation} 
which converges to zero as $n_1n_2 \rightarrow \infty$ by assumption and the dominated convergence theorem. Hence it remains to show that 
 \begin{equation} \label{tightaimboundedpart}
\lim_{ \lambda \to \infty} \limsup_{\min (n_1, n_2)  \to \infty}  \lambda^2   \p \Big (\max_{i \leq n_1 \atop j \leq n_2} |S'_{i,j}(f)|  \geq 5 \lambda \sqrt{n_1 n_2 } \Big )= 0 \, . 
\end{equation} 
With this aim, we start by noticing  that for any positive integer $q$ and $n$ such that $q \leq n$, setting $k_n = [n/q]$, we have 
\begin{equation} \label{tightp1}
\max_{ 1 \leq i \leq n} |s_i| \leq  \max_{ 1 \leq i \leq k_n} |s_{iq}| +   \max_{ 1 \leq j \leq k_n} \max_{ jq +1 \leq i \leq (j+1) q } |s_i  - s_{jq}|  \, .
\end{equation} 
Without loss of generality we can assume that $n_2 \leq n_1$.  Let  $q$ be a positive integer less than  $n_2$. Taking into account \eqref{tightp1}, we first write the following decomposition: 
\begin{multline} \label{tightp2}
\max_{i \leq n_1 \atop j \leq n_2} |S'_{i,j}(f)| \leq   \max_{i \leq k_{n_1} \atop j \leq n_2} |S'_{iq,j}(f)| + 
 \max_{0 \leq k \leq k_{n_1} \atop j \leq n_2}  \max_{kq+1 \leq i \leq (k+1) q}|S'_{i,j} (f)- S'_{kq,j}(f)| \\
\leq    \max_{ i \leq k_{n_1} \atop j \leq k_{n_2}} |S'_{iq,jq}(f)| +    \max_{i \leq k_{n_1} } \max_{0 \leq \ell \leq k_{n_2}}   \max_{\ell q+1 \leq j \leq (\ell+1) q}|S'_{iq,j} (f)- S'_{iq,\ell q} (f)|  \\
+   \max_{j \leq k_{n_2} } \max_{0 \leq k \leq k_{n_1}}   \max_{k q+1 \leq i \leq (k+1) q}|S'_{i,j q}(f) - S'_{kq,j  q}(f)|   +R_{n_1,n_2} (f)  \, .
\end{multline}
where 
\[
R_{n_1,n_2} (f) :=  \max_{0 \leq k \leq k_{n_1} \atop 
0 \leq \ell \leq k_{n_2}}  \max_{kq+1 \leq i \leq (k+1) q \atop \ell q+1 \leq j \leq (\ell+1) q } \Big | \sum_{u=kq +1}^i   \sum_{v=\ell q +1}^j X'_{u,v}   \Big | \, .
\]
It follows that 
\[
R_{n_1,n_2} (f) \leq  q^2 M \, .
\]
In what follows we shall select $q$  such that $ q^2 M \leq \lambda \sqrt{n_1 n_2}$ and to soothe the notation we write $S'_{i,j}$ for $S'_{i,j}(f)$.   Combined with \eqref{tightp2}, the selection of $q$ and $M$ will imply that 
\begin{equation} \label{tightp3}
\p \Big ( \max_{i \leq n_1 \atop j \leq n_2} |S'_{i,j}| \geq  5 \lambda   \sqrt{n_1n_2} \Big )  \leq   I_1 + I_2+I_3  \, .
\end{equation}
where 
\begin{align*} 
I_1 &:=  \p \Big (   \max_{ i \leq k_{n_1} \atop j \leq k_{n_2}} |S'_{iq,jq}|  \geq  2  \lambda \sqrt{n_1n_2} \Big )   \, , \\ I_2& :=   \p \Big (   \max_{i \leq k_{n_1} } \max_{0 \leq \ell \leq k_{n_2}}   \max_{\ell q+1 \leq j \leq (\ell+1) q}|S'_{iq,j} - S'_{iq,\ell q}|  \geq  \lambda \sqrt{n_1n_2} \Big )  \, ,  \\
 I_3& :=   \p \Big (    \max_{j \leq k_{n_2} } \max_{0 \leq k \leq k_{n_1}}   \max_{k q+1 \leq i \leq (k+1) q}|S'_{i,j q} - S'_{kq,j  q}|   \geq  \lambda \sqrt{n_1n_2} \Big )  \, .
\end{align*}
In the following, we select $q$ as follows
\[
q= \theta^{-1} \circ G^{-1} (v) \, .
\]
With this selection, taking into account \eqref{Rrightcontinuous}, it follows that 
\[
q^2M \leq   R(v) = R( R^{-1} (x) ) \leq x \, , 
\]
as desired since we can assume that $\lambda \geq 1$.   


\medskip

 In what follows we shall first assume that $q \leq  n_2$.  We first handle   the second term in the right-hand side of  \eqref{tightp3}.  With this aim,  we  notice that, by stationarity, 
\[
I_2 \leq ( k_{n_2}+1)    \p \Big (   \max_{1 \leq i \leq k_{n_1} }   \max_{1 \leq j \leq  q}|S'_{iq,j} |  \geq  \lambda \sqrt{n_1n_2} \Big )  \, .
\]
Hence, setting 
\[
U_{\ell} (j) = S'_{ \ell q, j} - S'_{ (\ell -1)q, j}  \, , 
\]
we derive that 
\begin{align} \label{dec1I2}
I_2 & \leq ( k_{n_2}+1)    \p \Big (   \max_{1 \leq i \leq [k_{n_1} /2 ]}   \max_{1 \leq j \leq  q} \Big | \sum_{\ell =1}^i U_{2 \ell} (j)  \Big |  \geq  \lambda \sqrt{n_1n_2}  /2\Big ) \nonumber   \\
& \quad \quad +   ( k_{n_2}+1)    \p \Big (   \max_{1 \leq i \leq[ (k_{n_1}-1) /2] }   \max_{1 \leq j \leq  q} \Big | \sum_{\ell =1}^i U_{2\ell +1} (j)  \Big |  \geq  \lambda \sqrt{n_1n_2}  /2\Big )  \nonumber  \\
& \leq I_{2,1} +  I_{2,2}  \, .
\end{align}
Let  
\[
{\tilde U}_{\ell} (j) = U_{2\ell} (j)  - \E ( U_{2\ell} (j)   | {\mathcal F}_{2(\ell-1) q,q}  ) \, ,
\]
where ${\mathcal F}_{2(\ell-1) q,q} = \sigma ( X_{u,v}, u \leq  2(\ell-1) q, v \leq q)$.
It follows that for any fixed integer $j$  in  $[1,q]$, $( {\tilde U}_{\ell} (j) )_{\ell \geq 1}$ is a stationary sequence of martingale differences with respect to $( {\mathcal G}_{\ell} )_{\ell \geq 1}$ where 
$ {\mathcal G}_{\ell}  = {\mathcal F}_{2\ell q,q} $.  Note first that 
\begin{multline}  \label{dec1I21-1}
 \p \Big (   \max_{1 \leq i \leq [ k_{n_1} /2] }   \max_{1 \leq j \leq  q} \Big | \sum_{\ell =1}^i  \E ( U_{2\ell} (j)   | {\mathcal F}_{2(\ell-1) q,q} )\Big |  \geq  \lambda \sqrt{n_1n_2}  /4\Big )  \\
 \leq \frac{4}{\lambda x }   \sum_{\ell =1}^{  [ k_{n_1} /2] } \sum_{u=(\ell -1) q +1}^{\ell q } \sum_{v=1}^q \Vert \E (X'_{u,v} |  {\mathcal F}_{ u-q,q}  ) \Vert_1 \ll  \frac{n_1 q }{\lambda x } \theta(q)   \, ,
\end{multline}
where we recall that $x= \sqrt{n_1n_2}$. 
Next, note that $  \big ( \max_{1 \leq j \leq q }  \big | \sum_{\ell =1}^i {\tilde U}_{\ell} (j)  \big |\big )_{i \geq 1}  $  is a submartingale with respect to $ (  {\mathcal G}_{i}  )_{i \geq 1}$.  Therefore, by Doob's maximal inequality, for any $p \geq 1$, 
\begin{equation}  \label{dec1I21-2}
 \p \Big (   \max_{1 \leq i \leq [ k_{n_1} /2] }   \max_{1 \leq j \leq  q} \Big | \sum_{\ell =1}^i  {\tilde U}_{\ell} (j)  \Big |  \geq  \lambda \sqrt{n_1n_2}  /4\Big )  \\
 \leq \frac{4^p}{\lambda^p x^p}   \Big  \Vert   \max_{1 \leq j \leq q }  \Big | \sum_{\ell =1}^{[k_{n_1} /2]} {\tilde U}_{\ell} (j)  \Big |\ \Big \Vert_p^p  \, .
\end{equation}
Let us choose $p=3$. Set $V_{\ell,t}:= \sum_{s=(2\ell -1)q +1}^{2 \ell q}  (  X'_{s,t}  - \E_{2(\ell-1) q,q} ( X'_{s,t}  )   )$. Note first that 
\[
 \max_{1 \leq j \leq q }  \Big | \sum_{\ell =1}^{[k_{n_1} /2]} {\tilde U}_{\ell} (j)  \Big | =  \max_{1 \leq j \leq q }  \Big | \sum_{t=1}^j   \sum_{\ell =1}^{[k_{n_1} /2]} V_{\ell,t} \Big | \, .
\]
Hence, by inequality  \eqref{ineRosRF} of the appendix and stationarity, we get 
\begin{multline}  \label{dec1I21-3}
\Big  \Vert   \max_{1 \leq j \leq q }  \Big | \sum_{\ell =1}^{[k_{n_1} /2]} {\tilde U}_{\ell} (j)  \Big |\ \Big \Vert_3^3   \ll  k_{n_1} q   \Vert S'_{q,1} \Vert_3^3  \\
+ k_{n_1} q  \Big (  \sum_{k=1}^{[k_{n_1} /2]}   \sum_{j=1}^{q } \frac{1}{k^{4/3}  j^{4/3}} 
 \Big \Vert   \sum_{\ell =1}^k \E_{2q,1}  ( {\tilde U}^2_{\ell} (j)   )\Big \Vert^{1/2}_{3/2} \Big )^{3} \, .
 \end{multline}
Now
\[
\E_{2q,1}  ( {\tilde U}^2_{\ell} (j)   )  \leq  \E_{2q,1}  ( U^{2}_{ 2 \ell} (j)   )   \leq   |  \E_{2q,1}  ( U^{2}_{ 2 \ell} (j)   )  - \E( U^{2}_{ 2 \ell} (j)   )  | +  \E( U^{2}_{ 2 \ell} (j)   )   \, .
\]
Next 
\[
 \sum_{\ell =1}^k  \E( U^{2}_{ 2 \ell} (j)  )   = k   \E( S^{\prime 2}_{q,j} ) \leq   k  \sum_{u,a=1}^q \sum_{v,b=1}^j  |\E ( X'_{u,v} X'_{a,b} ) |  \, .
\]
But, by \cite[Proposition 1]{DD03}, 
\begin{equation*} \label{directboundavectheta}
 |\E ( X'_{u,v} X'_{a,b} ) | \ll  \int_0^{\min (\theta (|a-u|),  \theta (|b-v|))} Q \circ G(u) du \, .
\end{equation*}
Therefore
\[
 \sum_{\ell =1}^k  \E( U^{2}_{ 2 \ell} (j)  )  \ll    k  q j  \int_0^1 R(u) Q(u) du  \, .
\]
On another hand, 
\begin{multline*} 
\Vert  \E_{1,1}  ( U^{2}_{ 2 \ell} (j)   )  - \E(U^{2}_{ 2 \ell} (j)   )   \Vert_{3/2}  \leq   \sum_{a,b=(2\ell -1)q +1} ^{2\ell q}\sum_{c,d=1}^j  \Vert  \E_{2q,1}  ( X'_{a,c}  X'_{b,d}  )  - \E(( X'_{a,c}  X'_{b,d}  )   \Vert_{3/2} \, .
\end{multline*}
But, 
for $\ell  \geq 2$  and  $j \leq q$, by stationarity and Lemma \ref{boundmomentp}, 
\begin{multline*}   \sum_{a,b=(2\ell -1)q +1} ^{2\ell q}\sum_{c,d=1}^j  \Vert  \E_{2q,1}  ( X'_{a,c}  X'_{b,d}  )  - \E(( X'_{a,c}  X'_{b,d}  )   \Vert_{3/2} \\
\leq  2^{2/3} q^2 j^2 M^{2/3}   \Big ( \int_0^{\theta((2\ell -2)q +1)} Q \circ G(u) du \Big )^{2/3}   \, .
\end{multline*}
Hence
\[
 \sum_{\ell =2}^k   \Vert  \E_{1,1}  ( U^{2}_{ 2 \ell} (j)   )  - \E( U^{2}_{ 2 \ell} (j)   )   \Vert_{3/2}   \ll  M^{2/3}   q^2 j^2  \sum_{\ell =2}^k  \Big ( \int_0^{\theta((2\ell -2)q +1)} Q \circ G(u) du \Big )^{2/3}  \, .
\]
But $\int_0^{\theta(k)} Q \circ G(u) du \leq C k^{-2}$. Hence $\sum_{\ell =2}^k \Big ( \int_0^{\theta((2\ell -2)q +1)} Q \circ G(u) du \Big )^{2/3}   \ll q^{-4/3}$. 
 It follows that 
 \[
   \sum_{\ell =2}^k  \Big \Vert   \E_{2q,1}  ( {\tilde U}^2_{\ell} (j)   )\Big \Vert_{3/2}  \ll M^{2/3}    j^2 q^{2/3}  + k q j  \, ,
 \]
 implying that 
\begin{equation}  \label{dec1I21-3*1}
  \sum_{k=1}^{[k_{n_1} /2]}   \sum_{j=1}^{q } \frac{1}{k^{4/3}  j^{4/3}} 
 \Big \Vert   \sum_{\ell =2}^k \E_{2q,1}  ( {\tilde U}^2_{\ell} (j)   )\Big \Vert^{1/2}_{3/2} \ll M^{1/3}    q +  q^{2/3} k^{1/6}_{n_1}  \, .
\end{equation} 
Next
\begin{equation}  \label{dec1I21-3*1prime}
\Vert   \E_{2q,1}  ( {\tilde U}^2_{1} (j)   ) \Vert_{3/2}   \leq   \Vert   \E_{2q,1}  ( U^{2}_{ 2 } (j)   )   \Vert_{3/2}  =   \Vert   \E_{\infty,1}  (  S^{\prime 2}_{q,j}  )   \Vert_{3/2}  \, .
\end{equation} 
 With the same arguments as  to prove Proposition \ref{proprecurrence}, we get 
 \[
  \Vert   \E_{\infty,1}  (  S^{\prime 2}_{q,j}  )   \Vert_{3/2} \ll  q^{2/3}  \Vert   \E_{1,1}  (  S^{\prime 2}_{1,j}  )   \Vert_{3/2} +  q^{2/3} \Big (  \sum_{k=0}^{r-1} 2^{-k  / 3}  \Vert  \E_{0,1} (   S^{\prime 2}_{2^k,j }   )  \Vert_{3/2}^{1/2}  \Big )^{2}   \, ,
\]
where $r$ is the unique integer such that $2^{r-1} \leq q < 2^r$. 
By stationarity, for $a \leq b$,  $ \Vert  \E_{0,1} (   S^{\prime 2}_{ a+b,j }   )  \Vert_{3/2} \leq  2   \Vert  \E_{0,1} (   S^{\prime 2}_{ a,j }   )  \Vert_{3/2}  +  2 \Vert  \E_{0,1} (   S^{\prime 2}_{ b-a,j }   )  \Vert_{3/2} $. So, by  Lemma 3.24 in \cite{MPU19},
\begin{equation}  \label{dec1I21-3*2}
  \Vert   \E_{\infty,1}  (  S^{\prime 2}_{q,j}  )   \Vert_{3/2} \ll  q^{2/3}  \Vert   \E_{1,1}  (  S^{\prime 2}_{1,j}  )   \Vert_{3/2} +  q^{2/3} \Big (  \sum_{i=1}^{q} i^{-4/3}  
   \Vert  \E_{0,1} (   S^{\prime 2}_{i,j }   )  \Vert_{3/2}^{1/2}  \Big )^{2}   \, . 
\end{equation}
Let us handle the quantity $\Vert  \E_{0,1} (   S^{\prime 2}_{k,\ell }   )  \Vert_{3/2}$. We first write
\[
\Vert  \E_{0,1} (   S^{\prime 2}_{k,\ell }   )  \Vert_{3/2} \leq \Vert \E_{ {\bf 1} } (  S^{\prime 2}_{{k}, \ell} ) - \E (  S^{\prime 2}_{{k}, \ell} )
\Vert_{3/2} + \E (  S^{\prime 2}_{{k}, \ell} )  \, .
\]
Since $\int_0^1 R(u) Q(u) du < \infty$, we infer that $\E (  S^{\prime 2}_{{k}, \ell} )  \ll   k \ell$. Next
\[
\Vert \E_{ {\bf 1} } (  S^{\prime 2}_{{k}, \ell} ) - \E (  S^{\prime 2}_{{k}, \ell} ) 
\Vert_{3/2}  \leq  \sum_{u,a=1}^k   \sum_{v,b=1}^\ell \Vert \E_{1,1} ( X'_{u,v} X'_{a,b} ) -   \E ( X'_{u,v} X'_{a,b} )  \Vert_{3/2}   \, .
\]
But, since the filtrations are commuting, 
\[
\Vert \E_{1,1} ( X'_{u,v} X'_{a,b} ) -   \E ( X'_{u,v} X'_{a,b} )  \Vert_{3/2}   \leq 2   \Vert \E_{a\wedge u, b \wedge v} (X'_{u,v} )   \E_{a\wedge u, b \wedge v}  (X'_{a,b} )  \Vert_{3/2}
 \, .
\]
By the first part of Lemma \ref{boundmomentp}, this implies that 
\[
\Vert \E_{1,1} ( X'_{u,v} X'_{a,b} ) -   \E ( X'_{u,v} X'_{a,b} )  \Vert^{3/2}_{3/2}   \ll M   \int_0^{ \min (  \theta  ( |a-u|) , \theta (|b-v)  ) } Q \circ G(u) du  \, .
\]
On another hand, using again  Lemma \ref{boundmomentp},
\[
\Vert \E_{1,1} ( X'_{u,v} X'_{a,b} ) -   \E ( X'_{u,v} X'_{a,b} )   \Vert^{3/2}_{3/2}    \ll M  \int_0^{  \min  ( \theta  ( u \wedge a ) , \theta (  v \wedge b )  ) } Q \circ G(u) du \, . 
\] 
So, overall,
\[
\Vert \E_{ {\bf 1} } (  S^{ \prime 2}_{{k}, \ell} ) - \E (  S^{ \prime  2}_{{k}, \ell} ) 
\Vert_{3/2}  \ll M^{2/3}\sum_{u=1}^k  \sum_{a=u}^k  \sum_{v=1}^\ell  \sum_{b=v}^{\ell}   A(u,a,v,b)  \, ,
\]
where
\[
A(u,a,v,b)  := \Big (  \int_0^{  \min   ( \theta  ( |a-u|) , \theta (|b-v) | , \theta (u) , \theta  ( v) ) } Q \circ G(u) du  \Big )^{2/3}\, .
\]
Now write that 
\begin{multline*}
 \sum_{u=1}^k  \sum_{a=u}^k  \sum_{v=1}^\ell  \sum_{b=v}^{\ell}   A(u,a,v,b) \leq  \sum_{u=1}^k  \sum_{a=u}^{2u}  \sum_{v=1}^\ell  \sum_{b=v}^{2v}   A(u,a,v,b) +  \sum_{u=1}^k  \sum_{a=2u+1}^k  \sum_{v=1}^\ell  \sum_{b=v}^{2v}   A(u,a,v,b)   \\ +  \sum_{u=1}^k  \sum_{a=u}^{2u}  \sum_{v=1}^\ell  \sum_{b=2v}^{\ell}   A(u,a,v,b)  +  \sum_{u=1}^k  \sum_{a=2u+1}^k  \sum_{v=1}^\ell  \sum_{b=2v+1}^{\ell}   A(u,a,v,b)  \, .
\end{multline*}
We handle the first term in the right-hand side, the others being controlled using similar arguments. For $k \geq \ell$, we have 
\begin{multline*}
 \sum_{u=1}^k  \sum_{a=u}^{2u}  \sum_{v=1}^\ell  \sum_{b=v}^{2v}   A(u,a,v,b)   =   \sum_{u=1}^k  \sum_{a=u}^{2u}  \sum_{v=1}^\ell  \sum_{b=v}^{2v}    \Big (  \int_0^{  \min  ( \theta  (u) , \theta  ( v) ) } Q \circ G(u) du  \Big )^{2/3}  \\
 \leq 2  \sum_{u=1}^\ell  \sum_{v=u}^\ell  u v    \Big (  \int_0^{   \theta  ( v)}  Q \circ G(u) du  \Big )^{2/3} + \ell^2    \sum_{u=\ell}^k  u    \Big (  \int_0^{   \theta  ( u)}  Q \circ G(u) du  \Big )^{2/3}   \ll \ell^2 k^{2/3}  \, ,
\end{multline*}
where for the last line we have used the fact that $ \int_0^{   \theta  ( u)}  Q \circ G(u) du = O (u^{-2})$. So, overall, 
\begin{equation}  \label{dec1I21-3*3}
\Vert \E_{ {\bf 1} } (  S^{ \prime 2}_{{k}, \ell} )
\Vert_{3/2} \ll      M^{2/3}   ( k \wedge \ell)^2  ( k \vee \ell)^{2/3}  +  k \ell \, .
\end{equation}
Hence, considering  \eqref{dec1I21-3*1prime}, \eqref{dec1I21-3*2} and \eqref{dec1I21-3*3}, we get 
\begin{multline}  \label{dec1I21-3*4}
 \sum_{k=1}^{[k_{n_1} /2]}   \sum_{j=1}^{q } \frac{1}{k^{4/3}  j^{4/3}} 
 \Big \Vert   \E_{2q,1}  ( {\tilde U}^2_{1} (j)   )\Big \Vert^{1/2}_{3/2} \\  \ll    q^{1/3} \sum_{j=1}^{q }  \sum_{i=1}^q \frac{(i j )^{1/2}+ M^{1/3} ( i \wedge j) ( i \vee j)^{1/3}  }{  i^{4/3}  j^{4/3}}  \ll M^{1/3}  q  + q^{2/3} \, .
 \end{multline}
So,  taking into account  \eqref{dec1I21-3*1} and \eqref{dec1I21-3*4}, we get 
 \[
k_{n_1} q  \Big (  \sum_{k=1}^{[k_{n_1} /2]}   \sum_{j=1}^{q } \frac{1}{k^{4/3}  j^{4/3}} 
 \Big \Vert   \sum_{\ell =1}^k \E_{2q,1}  ( {\tilde U}^2_{\ell} (j)   )\Big \Vert^{1/2}_{3/2} \Big )^{3} \ll   q^{3} k^{3/2}_{n_1} 
  +   n_1 M q^3      \, .
\]
Since $q^2 M \leq x$, it follows that 
\begin{equation}  \label{dec1I21-3*5}
k_{n_1} q  \Big (  \sum_{k=1}^{[k_{n_1} /2]}   \sum_{j=1}^{q } \frac{1}{k^{4/3}  j^{4/3}} 
 \Big \Vert   \sum_{\ell =1}^k \E_{2q,1}  ( {\tilde U}^2_{\ell} (j)   )\Big \Vert^{1/2}_{3/2} \Big )^{3} \ll   (n_1  q )^{3/2} +  x  n_1  q     \, .
 \end{equation}
We handle now the quantity  $  \Vert S'_{q,1} \Vert_3^3 $ of the right-hand side of \eqref{dec1I21-3}. Using inequality \eqref{ineRosRF} of the appendix, we get that
\[
 \Vert S'_{q,1} \Vert_3^3  \ll  q \Big ( M^{1/3}  +    \sum_{k=1}^{q}   \frac{\Vert \E_{ {\bf 1} } ( S^{\prime 2}_{{k}, 1} )
\Vert^{1/2}_{3/2} }{k^{4/3}  }  \Big  )^3 \, .
\]
Taking into account \eqref{dec1I21-3*3}, it follows that 
\begin{equation}  \label{dec1I21-3*6}
 k_{n_1} q  \Vert S'_{q,1} \Vert_3^3   \ll    n_1 ( M q^2 + q^{3/2})  \ll      x  n_1       + n_1  q^{3/2}   \, .
 \end{equation}
 Starting from \eqref{dec1I21-3} and considering the upper bounds  \eqref{dec1I21-3*5} and \eqref{dec1I21-3*6}, it follows that 
 \begin{equation}  \label{dec1I21-6}
 \Big  \Vert   \max_{1 \leq j \leq q }  \Big | \sum_{\ell =1}^{[k_{n_1} /2]} {\tilde U}_{\ell} (j)  \Big |\ \Big \Vert_3^3  \ll  (n_1  q )^{3/2} + x  n_1 q   \, .
 \end{equation}
 Taking into account \eqref{dec1I21-1},  \eqref{dec1I21-2},  \eqref{dec1I21-6} and the fact that $q \leq n_2$ and $x = \sqrt{n_1 n_2}$, we get 
\[
  I_{2,1}\ll   \frac{ x }{ \lambda}  \theta(q) +  \frac{1 }{ \lambda^3 }   \, . 
\]
 Clearly a similar upper bound holds for $  I_{2,2}$. It follows that 
\[
  \lambda^2 I_{2}\ll  \lambda x  \theta(q)   +   \lambda^{-1}   \, .
\]
 But $\theta(q) \leq  G^{-1} (v)  =  \int_0^v Q(u) du$ and $v =R^{-1} (x)$. Hence $\theta(q) \leq \int_0^1 Q(u)   {\bf 1}_{u < R^{-1} (x)} du = \int_0^1 Q(u) {\bf 1}_{x < R(u)} du$. This implies that 
 \begin{equation} \label{antedec1I21-7}
 x  \theta(q) \leq \int_0^1   Q(u) R(u) {\bf 1}_{R(u) >x} du \, .
 \end{equation}
 So, overall, 
  \begin{equation}  \label{dec1I21-7}
  \lambda^2 I_{2}\ll  \lambda \int_0^1 Q(u)  R(u) {\bf 1}_{R(u) >x} du   +   \lambda^{-1}   \, .
 \end{equation}
 Since  $ \int_0^1  Q(u) R(u)  du < \infty$ by assumption and $x \to \infty$ as $ n_1 \to \infty$, by the Lebesgue dominated theorem the first term in 
  the right-hand side of \eqref{dec1I21-7} is going to zero, as $ n_1 \to \infty$. Therefore
  \begin{equation}  \label{dec1I21-8}
\lim_{ \lambda \to \infty} \limsup_{\min(n_1,n_2) \to \infty}    \lambda^2 I_{2} = 0  \, .
 \end{equation}
 Proceeding similarly for $I_3$, we get  
   \begin{equation}  \label{dec1I21-9}
\lim_{ \lambda \to \infty} \limsup_{\min(n_1,n_2) \to \infty}    \lambda^2 I_{3} = 0  \, .
 \end{equation}
 It remains to handle the term $I_1$ in \eqref{tightp3}.  We start by writing that 
 \[
 I_1 \leq  \p \Big (   \max_{ i \leq k_{n_1} \atop j \leq k_{n_2}} |S_{iq,jq}|  \geq    \lambda \sqrt{n_1n_2} \Big )  +  \p \Big (   \max_{ i \leq k_{n_1} \atop j \leq k_{n_2}} |S''_{iq,jq}|  \geq    \lambda \sqrt{n_1n_2} \Big )  := I_{1,1} + I_{1,2} \, .
 \]
 But proceeding as to get \eqref{tightaimtruemoinsboundedpart}, with $x = \sqrt{n_1n_2}$, it follows that 
\[
  I_{1,2} \leq   \frac{2   }{ \lambda }  
  \int_0^1 R(u) Q(u)  {\bf 1}_{R(u) \geq  x }  du \, ,
\]
which converges to zero as $n_1n_2 \rightarrow \infty$ by assumption and the dominated convergence theorem.  It remains to handle the term $I_{1,1}$. 
With this aim,  for $N$  a fixed positive integer, we first  recall the  ortho-martingale coboundary decomposition \eqref{co-martdecwithN}: 
\begin{equation*} 
  f = m_N + (I - U)g_{1,N} + (I - V)g_{2,N}  + (I-U)(I-V)g_{3,N} + h_N \, .
\end{equation*}
Lemma 3.1 in \cite{VW14} states that any strictly stationary orthomartingale (suitably normalized) with respect to a commuting filtration and whose increments are in ${\mathbb L}^2$ is uniformly integrable. As a consequence, since  $(U^iV^j m_N)$ is a stationary field of orthomartingale differences in $  {\mathbb L}^2(\p)$, for any $N$ fixed,  
\begin{equation} \label{tightaimP0}
\lim_{ \lambda \to \infty} \limsup_{\min(n_1,n_2) \to \infty}  \lambda^2 \p \Big ( \max_{i \leq k_{n_1} \atop j \leq k_{n_2}} |S_{iq,jq} (m_N)  | \geq \lambda \sqrt{n_1n_2} \Big ) = 0 \, . 
\end{equation} 
Next, by stationarity, we infer that 
\begin{multline*}
\lambda^2 \p \Big ( \max_{i \leq k_{n_1} \atop j \leq k_{n_2}} |S_{iq,jq} ((I-U)g_{1,N})  | \geq \lambda \sqrt{n_1n_2} \Big )  \\
\leq  \lambda^2  \sum_{i=1}^{n_1}\p \Big ( \max_{1 \leq  j \leq n_2} \Big |  \sum_{\ell =1}^j  ( U - U^{i+1} ) V^{\ell } g_{1,N}  \Big  | \geq \lambda \sqrt{n_1n_2} \Big )  \\
\leq  \frac{8}{n_2} \Big \Vert   \max_{1 \leq  j \leq n_2} \Big |  \sum_{\ell =1}^j   V^{\ell } g_{1,N}  \Big  | 
{\mathbf 1}_{\{  \max_{1 \leq  j \leq n_2}  |  \sum_{\ell =1}^j   V^{\ell } g_{1,N}    |  \geq 2^{-1} \lambda \sqrt{n_1n_2}   \}} \Big \Vert_2^2 \, .
\end{multline*}
Since $(V^jg_{1,N})_j$ is a  stationary martingale differences sequence in $  {\mathbb L}^2(\p)$, its associated normalized partial sums are uniformly integrable. Therefore
\begin{equation} \label{tightaimP0bis}
\limsup_{\min(n_1,n_2) \to \infty}  \p \Big ( \max_{i \leq k_{n_1} \atop j \leq k_{n_2}} |S_{iq,jq} ((I-U)g_{1,N})  | \geq \lambda \sqrt{n_1n_2} \Big )  = 0 \, .
\end{equation}
Similarly, using that $(U^ig_{2,N})_i$ is a  stationary martingale differences sequence in $  {\mathbb L}^2(\p)$, 
\begin{equation} \label{tightaimP0ter}
\limsup_{\min(n_1,n_2) \to \infty}  \p \Big ( \max_{i \leq k_{n_1} \atop j \leq k_{n_2}} |S_{iq,jq} ((I-V)g_{2,N})  | \geq \lambda \sqrt{n_1n_2} \Big )  = 0 \, .
\end{equation}
Combining \eqref{tightaimP0}, \eqref{tightaimP0bis} and \eqref{tightaimP0ter}, it follows that 
\begin{equation*} \label{tightaimP1}
\lim_{ \lambda \to \infty} \limsup_{\min(n_1,n_2) \to \infty}  \lambda^2 \p \Big ( \max_{i \leq k_{n_1} \atop j \leq k_{n_2}} |S_{iq,jq} (f)  - S_{iq,jq} (h_N)  | \geq \lambda \sqrt{n_1n_2} \Big ) = 0 \, . 
\end{equation*} 
Hence, the theorem will follow if one can prove that 
\begin{equation} \label{tightaimP2}
 \lim_{N \rightarrow \infty}\limsup_{ \lambda \to \infty} \limsup_{\min(n_1,n_2) \to \infty}  \lambda^2 \p \Big ( \max_{i \leq k_{n_1} \atop j \leq k_{n_2}}| S_{iq,jq} (h_N)  | \geq \lambda \sqrt{n_1n_2} \Big ) = 0 \, . 
\end{equation}
We shall prove \eqref{tightaimP2}  with $\E_{0,0} ( U^N  f ) $ replacing $h_N$ (the proof with $\E_{0,0} ( V^N  f )$ and after with $ \E_{0,0} ( U^N V^N f )$ being similar). To simplify the notation, let $X^{(N)}_{i,j} =U^i V^j \E_{0,0} ( U^N  f ) $ and $S^{(N)}_{a,b} = \sum_{i=1}^{a} \sum_{j=1}^{b} X^{(N)}_{i,j}$.

Using the notation 
 \[
 U^{(N)}_{u,v} =  \sum_{i=(u-1)q}^{uq} \sum_{j=(v-1)q}^{vq} X^{(N)}_{i,j} \, ,
 \]
 we first consider the following decomposition: 
 \begin{multline} \label{tightaimP3dec}
  \max_{ i \leq k_{n_1} \atop j \leq k_{n_2}}  |S^{(N)}_{iq,jq}|  \leq    \max_{ i \leq [k_{n_1}/2]  \atop j \leq [k_{n_2}/2] }  \Big | \sum_{k=1}^i  \sum_{\ell=j}^i U^{(N)}_{2k,2\ell}    \Big |  +  \max_{ i \leq [(k_{n_1}+1)/2]  \atop j \leq [k_{n_2}/2] }  \Big | \sum_{k=1}^i  \sum_{\ell=j}^i U^{(N)}_{2k+1,2\ell}    \Big |  \\
  +  \max_{ i \leq [k_{n_1}/2]  \atop j \leq [(k_{n_2}+1)/2] }  \Big | \sum_{k=1}^i  \sum_{\ell=j}^i U^{(N)}_{2k,2\ell+1}    \Big |  +  \max_{ i \leq [(k_{n_1}+1)/2]  \atop j \leq [(k_{n_2}+1)/2] }  \Big | \sum_{k=1}^i  \sum_{\ell=j}^i U^{(N)}_{2k+1,2\ell+1}    \Big |  \, .
 \end{multline}
 We shall prove that 
 \begin{equation} \label{tightaimP3}
 \lim_{N \rightarrow \infty}\limsup_{ \lambda \to \infty} \limsup_{\min(n_1,n_2) \to \infty}  \lambda^2 \p \Big (  \max_{ i \leq [k_{n_1}/2]  \atop j \leq [k_{n_2}/2] }  \Big | \sum_{k=1}^i  \sum_{\ell=j}^i U^{(N)}_{2k,2\ell}    \Big | \geq \lambda \sqrt{n_1n_2} \Big ) = 0 \, . 
\end{equation}
The other terms in the right-hand side of \eqref{tightaimP3dec} can be handled similarly. To prove \eqref{tightaimP3}, we first write an orthomartingale decomposition. We set 
\[
 V_{k, \ell} =  U^{(N)}_{2k,2\ell}  - \E ( U^{(N)}_{2k, 2 \ell} | {\mathcal F}_{2 (k-1) q ,2 \ell q} )  -  \E (  U^{(N)}_{2k,2\ell} | {\mathcal F}_{2 k q ,2 ( \ell -1) q} ) + \E ( U^{(N)}_{2k, 2 \ell} | {\mathcal F}_{2 (k-1) q ,2 (\ell -1) q} )   \, ,
 \]
 and we note that 
  \begin{multline*} 
 \p \Big (  \max_{ i \leq [k_{n_1}/2]  \atop j \leq [k_{n_2}/2] }  \Big | \sum_{k=1}^i  \sum_{\ell=j}^i  ( U^{(N)}_{2k,2\ell} -   V_{k, \ell} )  \Big | \geq \lambda \sqrt{n_1n_2} \Big )    \\
 \leq \frac{1}{\lambda \sqrt{n_1n_2} }  \sum_{k=1}^{ [k_{n_1}/2] } \sum_{\ell=1}^{ [k_{n_2}/2] }  \Big (  
 \Vert  \E ( U^{(N)}_{2k, 2 \ell} | {\mathcal F}_{2 (k-1) q ,2 \ell q} )  \Vert_1  \\ + \Vert  \E (  U^{(N)}_{2k,2\ell} | {\mathcal F}_{2 k q ,2 ( \ell -1) q} )  \Vert_1+  \Vert \E ( U^{(N)}_{2k, 2 \ell} | {\mathcal F}_{2 (k-1) q ,2 (\ell -1) q} )  \Vert_1 \Big )\, .
 \end{multline*}
 By stationarity and \eqref{antedec1I21-7}, it follows that 
  \begin{multline*} 
 \p \Big (  \max_{ i \leq [k_{n_1}/2]  \atop j \leq [k_{n_2}/2] }  \Big | \sum_{k=1}^i  \sum_{\ell=j}^i  ( U^{(N)}_{2k,2\ell} -   V_{k, \ell} )  \Big | \geq \lambda \sqrt{n_1n_2} \Big )    \\
 \leq \frac{3}{\lambda \sqrt{n_1n_2} } n_1 n_2  \theta(q) \leq  \frac{3 x}{\lambda  } \theta(q)  \leq  3 \lambda^{-1}\int_0^1  Q(u) R(u) {\bf 1}_{R(u) >x} du \, ,
 \end{multline*}
which converges to zero as $ n_1 \to \infty$ since  $ \int_0^1   Q(u) R(u)  du < \infty$ by assumption.  Hence to prove \eqref{tightaimP3}, it remains to show that 
\begin{equation} \label{tightaimP3bis}
 \lim_{N \rightarrow \infty}\limsup_{ \lambda \to \infty} \limsup_{\min(n_1,n_2) \to \infty}  \lambda^2 \p \Big (  \max_{ i \leq [k_{n_1}/2]  \atop j \leq [k_{n_2}/2] }  \Big | \sum_{k=1}^i  \sum_{\ell=j}^i V_{k,\ell}    \Big | \geq \lambda \sqrt{n_1n_2} \Big ) = 0 \, . 
\end{equation}
By \cite[Prop. 2.2.1]{K02}, it follows that the sequences $ \big (  \max_{ i \leq [k_{n_1}/2]  }
 \big | \sum_{k=1}^i  \sum_{\ell=j}^i V_{k,\ell}    \big |  \big )_{j \geq 1}$ and $ \big (  \max_{  j \leq [k_{n_2}/2] } 
 \big | \sum_{k=1}^i  \sum_{\ell=j}^i V_{k,\ell}    \big |  \big )_{i \geq 1}$ are both one parameter submartingales. Therefore, applying twice Doob's maximal inequality, we get
\[
 \p \Big (  \max_{ i \leq [k_{n_1}/2]  \atop j \leq [k_{n_2}/2] }  \Big | \sum_{k=1}^i  \sum_{\ell=j}^i V_{k,\ell}    \Big | \geq \lambda \sqrt{n_1n_2} \Big ) 
 \leq \frac{4}{\lambda^2 n_1n_2 }  \sum_{k=1}^{ [k_{n_1}/2] } \sum_{\ell=1}^{ [k_{n_2}/2] }  
 \Vert  V_{k,\ell}  \Vert_2^2 \, .
\]
By stationarity
\[
\Vert  V_{k,\ell}  \Vert_2^2  \leq  \Vert  U^{(N)}_{k,\ell}  \Vert_2^2 = \Vert S^{(N)}_{q,q} \Vert_2^2\, .
\]
Since $S^{(N)}_{q,q} = \sum_{i,j}^{q} X^{(N)}_{i,j}$, using the properties of commuting filtrations, it follows that 
\[
\Vert S^{(N)}_{q,q} \Vert_2^2 \leq 2 \sum_{i,j,u,v=1}^q   \Vert \E_{i \wedge u, j \wedge v}  ( X^{(N)}_{i,j}) \E_{i \wedge u, j \wedge v}  ( X^{(N)}_{u,v}) \Vert_1 \, .
\]
Hence by Proposition 1 in \cite{DD03} and inequality (4.6) in \cite{Rio17}
\[
\Vert S^{(N)}_{q,q} \Vert_2^2 \leq 2 q^2 \sum_{a=N}^q  \sum_{b=1}^q  \int_0^{ \theta(a) \wedge \theta(b) } Q \circ G(u) du   \leq 6  q^2 \sum_{a=N}^q  a \int_0^{ \theta(a) \ } Q \circ G(u) du \, .
\]
So, overall, 
\[
\lambda^2 \p \Big (  \max_{ i \leq [k_{n_1}/2]  \atop j \leq [k_{n_2}/2] }  \Big | \sum_{k=1}^i  \sum_{\ell=j}^i V_{k,\ell}    \Big | \geq \lambda \sqrt{n_1n_2} \Big ) 
 \leq  24  \sum_{a\geq N}   a \int_0^{ \theta(a) } Q \circ G(u) du \, ,
\]
which converges to zero as $N \rightarrow \infty$ by assumption, so  \eqref{tightaimP3bis} is proved.   The proof of  \eqref{tightaim} is then complete when $q \leq n_2$. 

\medskip

Assume now that $q > n_2$ (recall that $q < \sqrt{n_1n_2}$) and let us highlight the differences compared to the previous case. Instead of \eqref{tightp2}, we write 
\begin{equation} \label{tightp2bis}
\max_{i \leq n_1 \atop j \leq n_2} |S'_{i,j}| \leq    \max_{i \leq k_{n_1} \atop j \leq n_2} |S'_{iq,j}| +  2 q n_2 M   \, ,
\end{equation}
and note that  $qn_2M \leq q^2 M \leq x$.  In addition,
\begin{multline*}
 \p \Big (   \max_{i \leq k_{n_1} \atop j \leq n_2} |S'_{iq,j}|  > \lambda \sqrt{n_1n_2 }\Big )  \leq   \p \Big (   \max_{1 \leq i \leq [k_{n_1} /2 ]}   \max_{1 \leq j \leq  n_2} \Big | \sum_{\ell =1}^i U_{2 \ell} (j)  \Big |  \geq  \lambda \sqrt{n_1n_2}  /2\Big ) \\
+     \p \Big (   \max_{1 \leq i \leq[ (k_{n_1}-1) /2] }   \max_{1 \leq j \leq  n_2} \Big | \sum_{\ell =1}^i U_{2\ell +1} (j)  \Big |  \geq  \lambda \sqrt{n_1n_2}  /2\Big )  \, .
\end{multline*}
Proceeding as to handle $I_2$ in the first part of the proof, we get 
\begin{multline*}
\lambda^2  \p \Big (   \max_{i \leq k_{n_1} \atop j \leq n_2} |S'_{iq,j}|  > \lambda \sqrt{n_1n_2 }\Big )  \leq   \lambda   \int_0^1 Q(u) R(u)  {\bf 1}_{R(u) >x}  du   \\   +  \frac{k_{n_1} n_2}{\lambda (n_1n_2)^{3/2}} \Big ( M q n_2 + q^{3/2} n_2 ^{1/2} k_{n_1}^{1/2}  +n_2 q \Big )       \, .
\end{multline*}
Hence, recalling that $q >n_2$ and that $q^2 M \leq x = \sqrt{n_1n_2}$, we derive 
\[
\lambda^2  \p \Big (   \max_{i \leq k_{n_1} \atop j \leq n_2} |S'_{iq,j}|  > \lambda \sqrt{n_1n_2 }\Big )  \leq  \lambda   \int_0^1 Q(u) R(u)  {\bf 1}_{R(u) >x}  du + \lambda^{-1}  \Big ( 1 +  \frac{1}{q} + \sqrt{ \frac{n_2}{n_1}}  \Big )  \, .\]
Therefore
\[
 \lim_{ \lambda \to \infty} \limsup_{\min(n_1,n_2) \to \infty}  \p \Big (   \max_{i \leq k_{n_1} \atop j \leq n_2} |S'_{iq,j}|  > \lambda \sqrt{n_1n_2 }\Big ) =0  \, ,
\]
ending the proof of  \eqref{tightaim} when $q > n_2$. $\square$

\subsection{Proof of Theorem \ref{FCLTthmreverse}}

The convergence of the finite dimensional laws follows from the arguments developed in Step 1 of the proof of Theorem \ref{FCLTthm}  combined with those given in the proof of Theorem \ref{TCLCFreverse}. Hence Theorem \ref{FCLTthmreverse} will follow provided we can prove the tightness of the partial sums process, which when  $d=2$ reads as 
\begin{equation}
\label{tightaimReverse}
\lim_{ \lambda \to \infty} \limsup_{n_1 \wedge n_2 \to \infty}  \lambda^2 
  \p \Big ( \max_{k \leq n_1, \ell \leq n_2} \big  | \sum_{i=1}^k \sum_{j = 1}^{\ell}  U_1^{n_1-i+1} U_2^{n_2-j+1}  f  \big | \geq \lambda \sqrt{n_1 n_2 } \Big ) = 0 \, . 
\end{equation} 
This can be done by following the lines of the proof of Proposition \ref{Proptight},   replacing $X'_{i,j}$ by $U_1^{n_1-i+1} U_2^{n_2-j+1}  f $ and 
$\sigma ( X_{i,j} , i \leq u, j \leq v)$ by $  {\mathcal F}_{n_1-u+1,n_2-v+1}  = {\mathcal F}^{(1)}_{n_1-u+1}  \cap   {\mathcal F}^{(2)}_{n_2-v+1} $.  Hence, we infer that  \eqref{tightaimReverse} holds provided that \eqref{plusieurs-points} and 
\eqref{condplusieurs-points}  are satisfied.  $\square$


\section{Appendix}

\setcounter{equation}{0}

In this section,  $(X_{{\bf i}})_{{\bf i} \in {\mathbb Z}^d}$ is a strictly stationary and centered random field adapted to a stationary filtration $({\mathcal F}_{{\bf i}})_{{\bf i} \in {\mathbb Z}^d}$. For any multi-integer $ {\bf i} $ of $ {\mathbb Z}^d$, recall that the notation $ \E_{ {\bf i} } ( X)$ means $ \E ( X | {\mathcal F}_{{\bf i}})$.  Moreover, if $d$ is a positive integer, we denote by $[d]$ the set $\{1, \ldots, d\}$ and,  for two elements ${\bf i}=(i_q)_{1 \leq q \leq d}$ and ${\bf j}=(i_q)_{1 \leq q \leq d}$ of ${\mathbb Z}^d$, ${\bf i} \preceq {\bf j}$ means  $i_q \leq j_q$ for any $q \in [d]$.

The following Rosenthal-type inequality will be  useful to prove tightness. It is an extension to the random fields setting of Theorem 6 in Merlev\`ede-Peligrad \cite{MP13}. 
\begin{theorem} \label{RosenthalRF}
Let $p  >2$ and $\delta = \min  \big ( 1/2, 1/(p-2) \big )$.  There exists a positive constant $K_{p,d}$ depending on $(p,d)$, such that for any ${\bf n}=(n_1, \ldots, n_d)$, 
\begin{equation*}
\Vert  \max_{{\bf k}  \preceq {\bf n}} |S_{\bf k} | \Vert_p \leq K_{p,d} \prod_{i=1}^d 2^{r_i/p} \left \{  \Vert X_{{\bf 1}}\Vert_p  +   \left [    \sum_{k_1=0}^{r_1-1} \ldots   \sum_{k_d=0}^{r_d-1}  \frac{\Vert \E_{ {\bf 1} } ( S^2_{2^{k_1}, \ldots, 2^{k_d}} )
\Vert^{\delta}_{p/2} }{2^{2 \delta k_1/p}  \cdots 2^{ 2 \delta k_d/p}}   \right ]^{1/(2 \delta)}  \right \}  \, ,
\end{equation*}
where for any $i \in \{1, \ldots, d \}$, $r_i$ is the unique integer such that $2^{r_i-1} \leq n_i < 2^{r_i}$. 
\end{theorem} 
By stationarity, since \[ \Vert \E_{ {\bf 1} } ( S^2_{{k_1}, \ldots, n+m, \ldots {k_d}} )
\Vert^\delta_{p/2} \leq 2^\delta \Vert \E_{ {\bf 1} } ( S^2_{{k_1}, \ldots, n, \ldots {k_d}} )
\Vert_{p/2}^\delta + 2^\delta \Vert \E_{ {\bf 1} } ( S^2_{{k_1}, \ldots, m, \ldots {k_d}} )
\Vert_{p/2}^\delta \, , \]applying Lemma 3.24 in \cite{MPU19} several times, we also have 
\begin{equation} \label{ineRosRF}
\Vert  \max_{{\bf k}  \preceq {\bf n}} |S_{\bf k} | \Vert_p \leq K'_{p,d} \prod_{i=1}^dn_i^{1/p} \left \{  \Vert X_{{\bf 1}}\Vert_p  +   \left [    \sum_{k_1=1}^{n_1} \ldots   \sum_{k_d=1}^{n_d}  \frac{\Vert \E_{ {\bf 1} } ( S^2_{{k_1}, \ldots, {k_d}} )
\Vert^{\delta}_{p/2} }{k_1^{1+2 \delta/p}  \cdots { k^{1+2 \delta/p}_d}} \right ]^{1/(2 \delta)}  \right \} \, .
\end{equation}

A first step in the proof of Theorem \ref{RosenthalRF} will be to show that it is enough to prove the maximal inequality for $\max_{{\bf k}  \preceq {\bf n}} \Vert S_{\bf k}\Vert_p$. This is achieved by using Proposition \ref{extensionmaxdyadic}  below. This proposition is an extension to stationary random fields of the maximal inequality (7) in \cite{MP13} (see also Proposition 3.14 and its corollary in \cite{MPU19}).  To state it, we need to introduce the following additional notations that are similar to those used in Giraudo \cite{Gi}.

$\bullet$ If $q \in [d]$, then ${\bf e}_q$ is the element of ${\mathbb N}^d$ such that the $q$th coordinate is equal to $1$, and all the
others to 0.

$\bullet$ For an element ${\bf k}=(k_q)_{1 \leq q \leq d}$ of ${\mathbb Z}^d$ and a non-empty subset $J$ of $ [d]$, the multiindexes ${\bf k}_J $  and 
${\bf 2}^{{\bf k}_J }$ of ${\mathbb Z}^d $ are defined by ${\bf k}_J=\sum_{q \in J} k_q {\bf e}_q$ and ${\bf 2}^{{\bf k}_J}=\sum_{q \in J} 2^{k_q} {\bf e}_q$. Moreover $ | { \bf k}_J | := \prod_{\ell \in J} k_{\ell} $,  $|{\bf 2}^{{\bf k}_J }|:= \prod_{\ell \in J} 2^{k_{\ell}} $ and the notation $ {\bf k}_J \in [{\bf 1}_J , {\bf n}_J]  $ is used to mean that $ 1 \leq k_{\ell} \leq n_{\ell}$  for any $q \in J$.  For  $J$ a subset of $ [d]$,  we shall also denote $J^c$ the set $ [d] \backslash J$.

$\bullet$ For a filtration $({\mathcal F}_{{\bf i } })_{ {\bf i } \in {\mathbb Z}^d}$ and a subset $J$ of $ [d]$, we denote by ${\mathcal F}_{\infty {\bf 1}_J}$ the $\sigma$-algebra generated by 
$ \displaystyle \cup_{{\bf j} \in {\mathbb Z}^d, {\bf j}_{[d] \backslash J }  \preceq {\bf 0} } {\mathcal F}_{{\bf j } } $.  

\begin{proposition} \label{extensionmaxdyadic}
Let $p  \geq 2$ and $q=p/(p-1)$.  There exists a positive constant $C_{p,d}$ depending on $p$ and $d$, such that for any ${\bf r}=(r_1, \ldots, r_d)$, 
\begin{equation*}
\Vert  \max_{{\bf k}  \preceq {\bf 2}^{\bf r}}  | S_{\bf k} | \Vert_p  \leq  q^d  \Vert S_{{\bf 2}^{\bf r}}\Vert_p 
+ q^d  \sum_{J \subset [d] }  \prod_{i \in J} 2^{r_i/p}  \sum_{{\bf k}_J \in [{\bf 0}_J , {\bf r}_J -{\bf 1}_J] }     \frac{ \Vert \E( S_{ {\bf 2}^{{\bf r}_{ J^c} + { \bf k}_J} }|  {\mathcal F}_{\infty {\bf 1}_J} )
\Vert_{p} }{ | {\bf 2}^{{ \bf k}_J } |^{1/p} }    \, .
\end{equation*}
\end{proposition} 
We refer to the next inequality \eqref{prop18d=2} for the rewriting of the above inequality in case $d=2$. 

\medskip

Using stationarity, the following subadditivity property holds:
\[
\Vert \E( S_{  {\bf i }_{ J^c} + { \bf k}_J + { \bf \ell}_J} |  {\mathcal F}_{\infty {\bf 1}_J} )
\Vert_{p}  \leq \Vert \E( S_{  {\bf i }_{ J^c} + { \bf k}_J}  |  {\mathcal F}_{\infty {\bf 1}_J} )
\Vert_{p}  + \Vert \E( S_{  {\bf i }_{ J^c} +  { \bf \ell}_J} |  {\mathcal F}_{\infty {\bf 1}_J} )
\Vert_{p} \, .
\]
Hence, taking into account Lemma 3.24 in \cite{MPU19} several times, we get

\begin{corollary} \label{extensionmax}
Let $p  \geq 2$ and $q=p/(p-1)$.  There exists a positive constant $C_{p,d}$ depending on $p$ and $d$, such that for any ${\bf n}=(n_1, \ldots, n_d)$, 
\begin{multline*}
\Vert  \max_{{\bf k}  \preceq {\bf n}}  | S_{\bf k} | \Vert_p  \leq (2q)^d  \max_{{\bf k}  \preceq {\bf n}} \Vert S_{\bf k}\Vert_p  \\
+ C_{p,d}  \sum_{J \subset [d] }  \prod_{i \in J} n_i^{1/p}  \sum_{{\bf k}_J \in [{\bf 1}_J , {\bf n}_J] }  \max_{{\bf i}_{ J^c  } \in  [ 
{\bf 1}_{J^c  }, {\bf n}_{J^c  }] }     \frac{ \Vert \E( S_{ {\bf i}_{ J^c} + { \bf k}_J} |  {\mathcal F}_{\infty {\bf 1}_J} )
\Vert_{p} }{ | { \bf k}_J |^{1+1/p} }    \, .
\end{multline*} 
\end{corollary}

We refer to the next inequality \eqref{ineincased=2} for the rewriting of the above inequality in case $d=2$.

%

We end this section by proving Theorem \ref{RosenthalRF} and Proposition \ref{extensionmaxdyadic}.

\medskip

\noindent {\bf Proof of Theorem \ref{RosenthalRF}.}  We shall first establish the inequality for $\max_{{\bf k}  \preceq {\bf n}} \Vert S_{\bf k}\Vert_p$.  For any integer $i$ in $[1,d]$, let $r_i$ be the unique integer such that $2^{r_i-1} \leq n_i < 2^{r_i}$. Using the binary expansion of each integer and stationarity, we infer that 
\[
\max_{{\bf k} \preceq {\bf n}}\Vert S_{{\bf k}} \Vert_p \leq \sum_{k_1=0}^{r_1-1}  \dots  \sum_{k_d=0}^{r_d-1}  \Vert S_{2^{k_1}, \dots, 2^{k_d}} \Vert_p \, .
\]
Hence, we shall first prove the inequality for  $\Vert S_{{\bf 2}^{\bf r}} \Vert_p$, namely: there exists a positive constant $ C_{p,d}$ such that 
\begin{equation} \label{ineRosdya}
\Vert S_{{\bf 2}^{\bf r}} \Vert_p \leq C_{p,d} \prod_{i=1}^d 2^{r_i/p} \left \{  \Vert X_{{\bf 1}}\Vert_p  +   \left [    \sum_{k_1=0}^{r_1-1} \ldots   \sum_{k_d=0}^{r_d-1}  \frac{\Vert \E_{ {\bf 1} } ( S^2_{2^{k_1}, \ldots, 2^{k_d}} )
\Vert^{\delta}_{p/2} }{2^{2 \delta k_1/p}  \cdots 2^{ 2 \delta k_d/p}}   \right ]^{1/(2 \delta)}  \right \}  \, .
\end{equation}
This can be achieved  by induction on $d$. Notice first that for $d=1$, inequality   \eqref{ineRosdya} is Theorem 6 in Merlev\`ede-Peligrad \cite{MP13} (see their inequality (17) and use the fact that $\Vert \E_{0,2^s } ( S_{2^{k},2^{s}} )
\Vert_p \leq \Vert \E_{0,2^s } ( S^2_{2^{k},2^{s}} )
\Vert_{p/2}^{1/2}$). Now assuming that the inequality is true for any dimension $d' \leq d-1$,  proving that it holds for $d'=d$ can be done as for proving that it holds for $d=2$ knowing that it is true for $d=1$. Since this asks less notations and it is easier to follow we give below the steps allowing to show that  inequality \eqref{ineRosdya} holds with $d=2$.

By  inequality (17) in Merlev\`ede-Peligrad \cite{MP13} with $\delta$ defined in Theorem \ref{RosenthalRF} (and again that  
$\Vert \E_{0,2^s } ( S_{2^{k},2^{s}} )
\Vert^2_p \leq \Vert \E_{0,2^s } ( S^2_{2^{k},2^{s}} )
\Vert_{p/2}$), there exists a positive constant $C_p$ such that for any integers $r,s \geq 0$, we have 
\[
 \Vert S_{2^{r},2^{s}} \Vert_p \leq C_p 2^{r/p} \Bigl \{  \Vert S_{1,2^{s}} \Vert_p +    \Big (   \sum_{k=0}^{r-1} \frac{\Vert \E_{0,2^s } ( S^2_{2^{k},2^{s}} )
\Vert^{\delta}_{p/2} }{2^{ 2\delta k/p }}  \Big )^{1/(2 \delta) }\Bigr \}  \, .
\]
Using again \cite[Inequality (17)]{MP13} on the first term in the right-hand side, we get 
\begin{multline} \label{Pr1Rosen}
 \Vert S_{2^{r},2^{s}} \Vert_p \leq C^2_p 2^{(r+s)/p} \Bigl \{  \Vert X_{1,1}\Vert_p  +  \Big (   \sum_{\ell=0}^{s-1} \frac{\Vert \E_{1,0 } ( S^2_{1,2^{\ell}} )
\Vert^{ \delta}_{p/2} }{2^{ 2 \delta \ell/p}}  \Big )^{1/(2 \delta) } \Bigr \}  \\ +   C_p  2^{r/p}   \Big (   \sum_{k=0}^{r-1} \frac{\Vert \E_{0,2^s } ( S^2_{2^{k},2^{s}} )
\Vert^{\delta}_{p/2} }{2^{ 2\delta k/p }}  \Big )^{1/(2 \delta) } \, .
\end{multline}
To take care of the last term in the right-hand side of \eqref{Pr1Rosen}, we shall use the following proposition.

\begin{proposition} \label{proprecurrence} 
Let $q >1$ and $\kappa = \min  \big ( \frac{1}{2},   \frac{1}{2(q-1)} \big ) $. For any positive integers $k$ and $n$, we have  
\begin{equation} \label{ine1proprecurrence}
\Vert \E_{0,\infty} ( S^2_{k,2^n} ) \Vert_q^q   \leq 2^n \Big ( 2 \Vert \E_{0,1} ( S^2_{k,1} ) \Vert_q^q  + c_q^{q/\kappa} \Big (  \sum_{\ell=0}^{n-1} 2^{-\ell \kappa / q }  \Vert  \E_{0,0} (   { S}^2_{k,2^{\ell}}   )  \Vert_q^{\kappa}  \Big )^{q/\kappa} \Big )
\end{equation}
and then 
\begin{equation} \label{ine2proprecurrence}
\max_{j \leq n}  \Vert \E_{0,n} ( S^2_{k,j} ) \Vert_q^q   \leq \frac{2^s}{  ( {  (2^{1/(2q)} -1 )^{2q}}} \Big ( 2 \Vert \E_{0,1} ( S^2_{k,1} ) \Vert_q^q  + c_q^{q/\kappa} \Big (  \sum_{\ell=0}^{s-1} 2^{-\ell \kappa / q }  \Vert  \E_{0,0} (   { S}^2_{k,2^{\ell}}   )  \Vert_q^{\kappa}  \Big )^{q/\kappa} \Big )  \, ,
\end{equation}
where $c_q   $ is given in Lemma \ref{lmarecurrence} below and  $s$ is the unique integer such that $2^{s-1} \leq n < 2^s$. 
\end{proposition}
Before proving the proposition, we end the proof of Theorem \ref{RosenthalRF} in case $d=2$.  We apply Proposition \ref{proprecurrence}   with $q=p/2$ to the  last term of the right-hand side of \eqref{Pr1Rosen} and note that $\kappa = \delta$. Therefore
\begin{multline*}
 \Vert S_{2^{r},2^{s}} \Vert_p \leq C^2_p 2^{(r+s)/p} \Bigl \{  \Vert X_{1,1}\Vert_p  +   \Big (   \sum_{\ell=0}^{s-1} \frac{\Vert \E_{1,0 } ( S^2_{1,2^{\ell}} )
\Vert^{ \delta}_{p/2} }{2^{ 2 \delta \ell/p}}  \Big )^{1/(2 \delta) } \Bigr \}  \\ +   C_p  2^{(r+s)/p}    \left (   2^{2 \delta /p}   \sum_{k=0}^{r-1} \frac{\Vert \E_{0,1 } ( S^2_{2^{k},1} )
\Vert^{\delta}_{p/2} }{2^{ 2\delta k/p}}  +  c_{p/2}   \sum_{k=0}^{r-1}  \sum_{\ell=0}^{s-1}  \frac{\Vert \E_{0,0 } ( S^2_{2^{k},2^{\ell}} )
\Vert^{\delta}_{p/2} }{2^{  2\delta k/p}  2^{  2\delta \ell/p}}  \right )^{1/(2 \delta) } \, ,
\end{multline*}
which proves inequality \eqref{ineRosdya} in case $d=2$. 

To complete the proof of Theorem \ref{RosenthalRF}, we need to take care of the moments of the maximum of partial sums. This is achieved by combining  inequality \eqref{ineRosdya}  with  Proposition \ref{extensionmaxdyadic} provided we can suitably bound the terms  $ \Vert \E( S_{ {\bf 2}^{{\bf r}_{ J^c} + { \bf k}_J} }|  {\mathcal F}_{\infty {\bf 1}_J} )
\Vert_{p} $ with the help of quantities such as  $ \Vert \E_{{\bf 1}}( S^2_{ {\bf 2}^{k }} )
\Vert_{p/2} $. Let us give the details in case $d=2$. First, we clearly have
\[
 \Vert \E ( S_{2^k,2^\ell} |  {\mathcal F}_{0,0})  \Vert_p \leq  \Vert \E ( S^2_{2^k,2^\ell} |  {\mathcal F}_{0,0})  \Vert^{1/2}_{p/2} \, .
\]
Next, according to the arguments developed in the proof of Proposition \ref{proprecurrence}, 
\begin{multline*}
 \Vert \E ( S_{2^r,2^\ell} |  {\mathcal F}_{\infty, 0})  \Vert_p  \leq  \Vert \E ( S^2_{2^r,2^\ell} |  {\mathcal F}_{ \infty, 0})  \Vert^{1/2}_{p/2} \\
 \ll    2^{r/p}   \Big (  \sum_{k=0}^{r-1}   2^{-2k\delta/p}   \Vert \E ( S^2_{2^k,2^\ell} |  {\mathcal F}_{1,0})  \Vert^{\delta}_{p/2}  \Big )^{1/(2 \delta) } \, .
 \end{multline*}
 Since $\delta \leq 1/2$, this implies that 
 \[
 \sum_{\ell =0}^{s-1}    2^{(s-\ell)/p}  \Vert \E ( S_{2^r,2^\ell} |  {\mathcal F}_{\infty,0})  \Vert_p  \ll 
 2^{(r+s)/p}   \Big (  \sum_{\ell =0}^{s-1}     \sum_{k=0}^{r-1} 2^{-2\ell\delta/p}     2^{-2k\delta/p}   \Vert \E ( S^2_{2^k,2^\ell} |  {\mathcal F}_{1,0})  \Vert^{\delta}_{p/2}  \Big )^{1/(2 \delta) }
 \]
 and similarly
 \[
 \sum_{k =0}^{r-1}    2^{(r-k)/p}  \Vert \E ( S_{2^k,2^s} |  {\mathcal F}_{0, \infty})  \Vert_p  \ll 
 2^{(r+s)/p}   \Big (     \sum_{k=0}^{r-1}   \sum_{\ell =0}^{s-1}  2^{-2\ell\delta/p}     2^{-2k\delta/p}   \Vert \E ( S^2_{2^k,2^\ell} |  {\mathcal F}_{0,1})  \Vert^{\delta}_{p/2}  \Big )^{1/(2 \delta) } \, .
 \]
Since the filtrations are commuting, all these considerations end the proof of the maximal inequality stated in Theorem  \ref{RosenthalRF} provided Proposition \ref{proprecurrence} holds. 
\smallskip

To end the proof,  it remains to prove Proposition \ref{proprecurrence}. 
\medskip

\noindent {\bf Proof of Proposition \ref{proprecurrence}.}  Inequality \eqref{ine1proprecurrence} follows from the recurrence formula given in Lemma below. Let us introduce the following notation: For $q \geq 1$, let 
\[
a^q_{k,n}  = \Vert \E_{0,n} ( S^2_{k,n} ) \Vert_q^q \, .
\]
\begin{lemma} \label{lmarecurrence}
For $q \in ]1,2]$, we have 
\begin{equation} \label{decompositionakn2}
a^q_{k,2n}  \leq 2 a^q_{k,n} +   c_q a_{k,n}^{q- 1/2}  \Vert  \E_{0,0} (   { S}^2_{k,n}   )  \Vert_q^{1/2}  \, ,
\end{equation}
where $c_q = 5 q + 2^q$. 
For $q>2$, we have 
\begin{equation} \label{decompositionakn2q>2}
a^q_{k,2n}  \leq 2 a^q_{k,n} + c_q a_{k,n}^{q-1/(2q-2)}   \Vert \E_{0,0} (  S^2_{k,n}   ) \Vert^{1/{(2q-2)}}_q  \, ,
\end{equation}
where $c_q = 4^q ( 2+2^{q+1} ) +2^q$. 
\end{lemma}
Before proving the lemma, let us end the proof of Proposition  \ref{proprecurrence}. Let $\kappa = \min  \big ( \frac{1}{2},   \frac{1}{2(q-1)} \big ) $. From inequalities \eqref{decompositionakn2} or \eqref{decompositionakn2q>2}, by recurrence on the first term, we obtain for any 
positive integer $s$,
\[
a^q_{k,2^s}  \leq 2^s  \Big ( a^q_{k,1} + 2^{-1} c_q \sum_{\ell=0}^{r-1} 2^{-\ell} a_{k,2^{\ell}}^{q-\kappa}  \Vert  \E_{0,0} (   { S}^2_{k,2^{\ell}}   )  \Vert_q^{\kappa} \Big ) \, .
\]
With the notation $B_s = \max_{0 \leq \ell \leq s}  (a^q_{k,2^{\ell}}/2^{\ell} )$, it follows that 
\[
B_s \leq a^q_{k,1} + 2^{-1} c_q B_s^{1-\kappa/q } \sum_{\ell=0}^{s-1} 2^{- \kappa \ell/ q }  \Vert  \E_{0,0} (   { S}^2_{k,2^{\ell}}   )  \Vert_q^{\kappa}  \, .
\]
Taking into account that either $B_s \leq 2 a^q_{k,1}$ or $B_s^{\kappa/ q} \leq   c_q \sum_{\ell=0}^{s-1} 2^{-\ell \kappa / q }  \Vert  \E_{0,0} (   { S}^2_{k,2^{\ell}}   )  \Vert_q^{\kappa} $, we derive that 
\begin{equation} \label{consrecurrence}
a^q_{k,2^s} \leq  2^s  \Big ( 2 a^q_{k,1} + c_q^{q /\kappa }  \Big (  \sum_{\ell=0}^{s-1} 2^{-\ell \kappa / q}  \Vert  \E_{0,0} (   { S}^2_{k,2^{\ell}}   )  \Vert_q^{\kappa}  \Big )^{q / \kappa} \Big )  \, ,
\end{equation}
proving inequality \eqref{ine1proprecurrence}. To prove inequality \eqref{ine2proprecurrence}, we proceed as page 87 in \cite{MPU19}. Let $j$ be a positive integer and let  $s_j$ the positive integer such that $2^{s_j} -1 \leq n < 2^{s_j}$, we write the binary expansion of $j$ as follows
\[
j = \sum_{\ell=0}^{s_j-1} 2^\ell b_\ell \text{ where }b_{s_j-1}=1 \text{ and }b_\ell \in \{0,1 \} \text{ for }\ell=0, \dots,s_j-2 \, .
\]
Note that 
\[
S_{k, j } =  \sum_{\ell=0}^{s_j-1} b_\ell T_{k, 2^{\ell} } \text{ where }  T_{k, 2^{\ell} } = S_{k, n_\ell } - S_{k, n_{\ell-1}} \, , \,  n_{\ell} = \sum_{i=0}^{\ell} 2^i b_i \text{ and }n_{-1}=0 \, . 
\]
Therefore
\[
\E_{0, \infty} (S^2_{k, j } ) =   \sum_{\ell, \ell'=0}^{s_j-1} b_\ell  b_{\ell'} \E_{0, \infty} (T_{k, 2^{\ell} }  T_{k, 2^{\ell'} } ) \leq   \sum_{\ell, \ell'=0}^{s_j-1} b_\ell  b_{\ell'} \E^{1/2}_{0, \infty} (T^2_{k, 2^{\ell} } ) \E^{1/2}_{0, \infty} ( T_{k, 2^{\ell'} } ) \, , 
\]
implying, by Cauchy-Schwarz's inequality, that 
\begin{multline*}
 \Vert \E_{0, \infty} (S^2_{k, j } ) \Vert_q  \leq   \sum_{\ell, \ell'=0}^{s_j-1}  \Vert \E^{1/2}_{0, \infty} (T^2_{k, 2^{\ell} } ) \E^{1/2}_{0, \infty} ( T_{k, 2^{\ell'} } )  \Vert_q 
 \\ \leq   \sum_{\ell, \ell'=0}^{s_j-1}  \Vert \E_{0, \infty} (T^2_{k, 2^{\ell} } ) \Vert_q^{1/2}  \Vert \E_{0, \infty} (T^2_{k, 2^{\ell'} } ) \Vert_q^{1/2} \, . 
\end{multline*}
Hence,  using stationarity,
\[
 \Vert \E_{0, \infty} (S^2_{k, j } ) \Vert_q  \leq   \Big (    \sum_{\ell=0}^{s_j-1}  \Vert \E_{0, \infty} (S^2_{k, 2^{\ell} } ) \Vert_q^{1/2}   \Big )^2 \, . 
\]
Inequality \eqref{ine2proprecurrence} follows by taking into account \eqref{ine1proprecurrence} in the inequality above together with the fact that 
\[
 \Big (   \sum_{\ell=0}^{s_j-1}  2^{\ell /(2q)} \Big )^{2q} =   \Big (   \frac{  2^{s_j  /(2q)}  - 1 } { 2^{1/(2q)} -1 }\Big )^{2q}  \leq   \frac{  2^{s_j }  } {  (2^{1/(2q)} -1 )^{2q}} \, .
\]
To end the proof of Proposition \ref{proprecurrence}, it remains to prove Lemma \ref{lmarecurrence}.

\smallskip

\noindent {\bf Proof of Lemma \ref{lmarecurrence}.} We start with the following inequality: For any $u,v \in {\mathbb R}^+$ and any $z \in {\mathbb R}$, note that 
\[
| u+v +z|^q\leq  | u+v |^q + q |z| ( u+v )^{q-1} + | z|^q \, .
\]
Indeed by the Taylor formula at order $1$ applied with the fonction $h(x) = |x|^q$ and since $h'(x) = q |x|^{q-1} sign(x)$, we have 
\[
| u+v +z|^q -   | u+v |^q  \leq  q  |z| (u+v)^{q-1} + |z| q  \int_0^1 \big ( ( u+v + t |z| )^{q-1} - ( u+v )^{q-1} \big ) dt \, .
\]
Since $q-1 \in [0,1]$,  $|  ( u+v + t |z|)^{q-1} - ( u+v )^{q-1}|  \leq  ( t |z|)^{q-1} $, it follows that
\[
| u+v +z|^q -   | u+v |^q  \leq  q  |z| (u+v)^{q-1} + |z|^q \, .
\]
Next, again by the Taylor formula at order $1$,
\[
( u+v )^q -    u^q  = q v u^{q-1} + q v  \int_0^1 \big ( ( u+ tv )^{q-1} -u^{q-1} \big ) dt  \leq q v u^{q-1} + v^q \, .
\]
So, overall, we get that for any $u,v \in {\mathbb R}^+$ and any $z \in {\mathbb R}$,
\begin{equation} \label{ineqdirect}
| u+v + z  |^q -    u^q -v^q  \leq q  |z| (u+v)^{q-1} + |z|^q + q v u^{q-1} \, .
\end{equation}

We go back to the proof of inequality \eqref{decompositionakn2}. Let ${\bar S}_{k,n} = S_{k,2n} - S_{k,n}$. Write first that
\[
\E_{0,2n} ( S^2_{k,2n} ) = \E_{0,2n} (  ( S_{k,n}  + {\bar S}_{k,n} )^2) = \E_{0,2n} (  S^2_{k,n}   ) +   \E_{0,2n} ( {\bar S}^2_{k,n} ) + 2   \E_{0,2n} ( S_{k,n}   {\bar S}_{k,n} )  \, .
\]
Since  $S^2_{k,n}$ is ${\mathcal F}_{k,n}$-measurable,  and the filtrations are commuting 
\[
\E_{0,2n} (  S^2_{k,n}   )  = \E_{0,2n} (   \E_{k,n}  (S^2_{k,n}  ) )  = \E_{0,n} (  S^2_{k,n}   )  \, .
\]
We apply now Inequality \eqref{ineqdirect} with $u =  \E_{0,n} (  S^2_{k,n}   )$, $v =   \E_{0,2n} ( {\bar S}^2_{k,n} )$ and $z = 2   \E_{0,2n} ( S_{k,n}   {\bar S}_{k,n} )$, and we take the expectation. We get 
\begin{equation} \label{decompositionakn}
a^q_{k,2n}  \leq \Vert  \E_{0,n} (  S^2_{k,n}   ) \Vert_q^q +  \Vert   \E_{0,2n} ( {\bar S}^2_{k,n} )\Vert_q^q +  2 q A + 2q B+  2^q C +  q D \, ,  
\end{equation}
where 
\begin{equation} \label{defAB}
A =   \E \Big (    \E_{0,2n} (  | S_{k,n}   {\bar S}_{k,n}  | ) 
(  \E_{0,n} (  S^2_{k,n}   )  )^{q-1}  \Big )  \, ,  \, B =   \E \Big (    \E_{0,2n} (  | S_{k,n}   {\bar S}_{k,n}  | ) 
(  \E_{0,2n} ( {\bar S}^2_{k,n} )  )^{q-1}  \Big )  \, , 
\end{equation}
\begin{equation} \label{defCD}
C =   \Vert    \E_{0,2n} ( S_{k,n}   {\bar S}_{k,n} )  \Vert_q^q \mbox{ and } D= \E \Big ( \E_{0,2n} ( {\bar S}^2_{k,n} )    ( \E_{0,n} (  S^2_{k,n}   ) )^{q-1}\Big )  \, .
\end{equation}
Let us analyze each of the terms $A$, $B$, $C$ and  $D$. We have 
\[
A =  \E \Big (    \E_{0,2n} (  | S_{k,n}   {\bar S}_{k,n}  | ) 
(  \E_{0,n} (  S^2_{k,n}   )  )^{q-1}  \Big )  = \E \Big (    \E_{0,n} (  | S_{k,n}   {\bar S}_{k,n}  | ) 
(  \E_{0,n} (  S^2_{k,n}   )  )^{q-1}  \Big )   \, .
\]
Hence, by H\"older's inequality twice and stationarity, 
\begin{multline} \label{borneA}
A 
\leq  \E \Big (    \E^{1/2}_{0,n} (   {\bar S}^2_{k,n}   ) 
(  \E_{0,n} (  S^2_{k,n}   )  )^{q-1/2}  \Big )   \\ \leq   \Vert   \E_{0,n} (   {\bar S}^2_{k,n}   )  \Vert_q^{1/2}
 \Vert \E_{0,n} (  S^2_{k,n}   )   \Vert_q^{q-1/2} = \Vert  \E_{0,0} (   { S}^2_{k,n}   )  \Vert_q^{1/2}
a_{k,n}^{q-1/2}  \, .
\end{multline}
To take care of $B$, we  set 
\[
\alpha=   \left (  \frac{ \Vert \E_{0,0} (  S^2_{k,n}  )  \Vert_q }{a_{k,n}} \right ) ^{q-3/2}  \text{ and } \eta = \left (  \frac{   \E_{0,n} (   {\bar S}^2_{k,n}   )}{  \E_{0,2n} (   {\bar S}^2_{k,n}   )} \right )^{2-q} \, ,
\]
with the convention that $0/0=1$. Clearly we can assume that  $\alpha >0$, otherwise there is nothing to prove. Since $\eta$ is ${{\mathcal F}}_{0,2n}$ measurable and $ \E_{0,2n} (  S^2_{k,n})=  \E_{0,n} (  S^2_{k,n})$,  we start by writing that 
\begin{multline*}
2 B \leq \alpha^{-1}   \E \Big ( \eta^{-1} \E_{0,n} (   S^2_{k,n}  ) 
(  \E_{0,2n} ( {\bar S}^2_{k,n} )  )^{q-1}  {\bf 1}_{\{ \E_{0,n} (   {\bar S}^2_{k,n}   ) \neq 0 \}}   \Big ) \\ +  \alpha   \E \Big ( \eta 
(  \E_{0,2n} ( {\bar S}^2_{k,n} )  )^{q}   {\bf 1}_{  \{\E_{0,2n} (   {\bar S}^2_{k,n}   ) \neq 0 \}}\Big )  =: \alpha^{-1} B_1 + \alpha  B_2 \, .
\end{multline*}
By the definition of $\eta$,  H\"older's inequality  and stationarity, we get 
\begin{multline*}
B_1 =  \E \Big (  \E_{0,n} (   S^2_{k,n}  ) \E^{q-2}_{0,n} (    {\bar S}^2_{k,n}  ) 
  \E_{0,2n} ( {\bar S}^2_{k,n} ) {\bf 1}_{ \{ \E_{0,n} (   {\bar S}^2_{k,n}   ) \neq 0 \}}  \Big )   \\
  = \E \Big (  \E_{0,n} (   S^2_{k,n}  ) \E^{q-2}_{0,n} (    {\bar S}^2_{k,n}  ) 
  \E_{0,n} ( {\bar S}^2_{k,n} )   {\bf 1}_{ \{ \E_{0,n} (   {\bar S}^2_{k,n}   ) \neq 0 \}}   \Big )  \\ =  \E \Big (  \E_{0,n} (   S^2_{k,n}  ) \E^{q-1}_{0,n} (    {\bar S}^2_{k,n}  ) \Big ) 
  \leq   a_{k,n}   \Vert  \E_{0,0} (  { S}^2_{k,n}  ) \Vert_q^{q-1} \, .
\end{multline*}
Next, by H\"older's inequality  and stationarity again, we derive 
\[
B_2 =  \E \Big (  \E^{2-q}_{0,n} (    {\bar S}^2_{k,n}  ) 
  \E^{2(q-1)}_{0,2n} ( {\bar S}^2_{k,n} )   \Big )  
  \leq   a^{2(q-1)}_{k,n}   \Vert  \E_{0,0} (  { S}^2_{k,n}  ) \Vert_q^{2-q} \, .
\]
So, overall, taking into account the definition of $\alpha$, it follows that 
  \begin{equation} \label{borneB}
  B \leq  a_{k,n}^{q-1/2} \Vert  \E_{0,0} (  { S}^2_{k,n}  ) \Vert_q^{1/2} 
 \, .
  \end{equation}

We handle now the term $C$ in \eqref{decompositionakn}.  We have 
\begin{multline} \label{forC}
C \leq  \E \Big (  \E_{0,2n} ( |S_{k,n}   {\bar S}_{k,n} |)  \E^{q-1}_{0,2n} (  | S_{k,n}   {\bar S}_{k,n}  | )  \Big )  \\
\leq 2^{-1}  \E \Big (  \E_{0,2n} ( |S_{k,n}   {\bar S}_{k,n} |)  \E^{q-1}_{0,n} (   S^2_{k,n}  )  \Big )  +  2^{-1}  \E \Big (  \E_{0,2n} ( |S_{k,n}   {\bar S}_{k,n} |)  \E^{q-1}_{0,2n} (   {\bar S}^2_{k,n}  )  \Big ) \\
\leq   2^{-1} (A+B) \, .
\end{multline}
Hence, taking into account \eqref{borneA}  and \eqref{borneB}, we get
 \begin{equation} \label{borneC}
  C \leq  a_{k,n}^{q-1/2}  \Vert  \E_{0,0} (   { S}^2_{k,n}   )  \Vert_q^{1/2}
 \, .
  \end{equation}
To handle the quantity $D$, we first notice that $D= \E \big ( \E_{0,n} ( {\bar S}^2_{k,n} )   \E^{q-1}_{0,n} (  S^2_{k,n}   )\big )$. Therefore, by H\"older's inequality, stationarity and the fact that 
$ \Vert \E_{0,0} (  { S}^2_{k,n}  ) \Vert_q \leq a_{k,n}$, we get 
 \begin{equation} \label{borneD}
D \leq  \Vert  \E_{0,0} (  { S}^2_{k,n}  ) \Vert_q a_{k,n}^{q-1}   \leq  \Vert  \E_{0,0} (  { S}^2_{k,n}  ) \Vert^{1/2}_q a_{k,n}^{q-1/2}   \, .
\end{equation}
Hence starting from \eqref{decompositionakn}, using stationarity  and taking into account the upper bounds \eqref{borneA}, \eqref{borneB}, \eqref{borneC}  and \eqref{borneD}, the inequality  \eqref{decompositionakn2} follows. 

\medskip

We turn now to the proof of \eqref {decompositionakn2q>2} so we consider the case $q >2$.  We shall apply twice inequality (87) in \cite{MP13}: for any $a,b \geq 0$, 
\[
(a+b)^q \leq a^q + b^q + 4^q ( a^{q-1} b + a b^{q-1}) \, , 
\]
first with $a= \E_{0,n} (  S^2_{k,n}   ) +   \E_{0,2n} ( {\bar S}^2_{k,n} ) $ and $b = 2  | \E_{0,2n} ( S_{k,n}   {\bar S}_{k,n} ) |$ and the second time with 
 $a= \E_{0,n} (  S^2_{k,n}   ) $ and $b=  \E_{0,2n} ( {\bar S}^2_{k,n} )$. Taking the expectation and using stationarity, we infer that 
 \begin{equation} \label{decompositionaknq>2}
a^q_{k,2n}  \leq 2 a^q_{k,n} +  4^q \times 2^{q-1}  ( A +  B) +  2^q C +  4^q D + 4^q \times 2^{q-1} (E+F)  + 4^q G \, ,  
\end{equation}
where $A,B,C, D$ are defined in \eqref{defAB} and \eqref{defCD},
\begin{equation} \label{defEF}
E =   \E \Big (    
(  \E_{0,n} (  S^2_{k,n}   )  )   \E^{q-1}_{0,2n} (  | S_{k,n}   {\bar S}_{k,n}  | )  \Big )  \, ,  \, F =  \E \Big (    
  \E_{0,2n} (   {\bar S}_{k,n}^2   )     \E^{q-1}_{0,2n} (  | S_{k,n}   {\bar S}_{k,n}  | )  \Big )   \, , 
\end{equation}
and
\begin{equation} \label{defG}
G = \E \big (   \E_{0,n} (  S^2_{k,n}   ) \E^{q-1}_{0,2n} ( {\bar S}^2_{k,n} ) \big )  \, .
\end{equation}
We start by giving an upper bound for $G$. Let $x = q-1$. Since $q>2$, $x >1$. By H\"older's inequality, we get
\begin{multline*}
G = \E \big (   \E_{0,n} (  S^2_{k,n}   )  \E^{1/x}_{0,2n} ( {\bar S}^2_{k,n} )  \E^{q-1-1/x}_{0,2n} ( {\bar S}^2_{k,n} ) \big )  \\
\leq  \E^{1/x} \Big [   \E^x_{0,n} (  S^2_{k,n}   )  \E_{0,2n} ( {\bar S}^2_{k,n} )  \Big ] \E^{(x-1)/x} \Big [    \E^{\frac{x(q-1)-1}{x-1}}_{0,2n} ( {\bar S}^2_{k,n} ) \Big ]  \, .  
\end{multline*}
By the property of commuting filtrations, H\"older's inequality and stationarity, 
\begin{align*}
 \E \big [   \E^x_{0,n} (  S^2_{k,n}   )  \E_{0,2n} ( {\bar S}^2_{k,n} )  \big ]  &=  \E \big [   \E^x_{0,n} (  S^2_{k,n}   )  \E_{0,n} ( {\bar S}^2_{k,n} )  \big ]  \\
& \leq \Vert \E_{0,n} (  S^2_{k,n}   ) \Vert^{q-1}_q \Vert \E_{0,0} (  S^2_{k,n}   ) \Vert_q  \, .
\end{align*}
Therefore, using that $\frac{x(q-1)-1}{x-1} = q$ and stationarity, we get 
\begin{equation} \label{borneG}
G\leq \Vert \E_{0,0} (  S^2_{k,n}   ) \Vert^{1/{(q-1)}}_q \Vert \E_{0,n} (  S^2_{k,n}   ) \Vert^{1+ q(q-2)/(q-1)}_q  = a_{k,n}^{q-1/(q-1)}   \Vert \E_{0,0} (  S^2_{k,n}   ) \Vert^{1/{(q-1)}}_q  \, .
\end{equation}
Next, for the quantity $A$ (resp. $D$), the estimate \eqref{borneA} (resp. \eqref{borneD}) is still available for $q>2$. Concerning $B$, we first write that for any $\gamma >0$, we have 
  \[
 2B \leq  \gamma^{-1}  \E \Big (    \E_{0,n} (  S^2_{k,n} )  (  \E_{0,2n} ( {\bar S}^2_{k,n} )  )^{q-1} \Big )  +  \gamma \E \Big (    
(  \E_{0,2n} ( {\bar S}^2_{k,n} )  )^{q}  \Big ) \leq   \gamma^{-1}  G +  \gamma a_{k,n}^q \, .
 \]
Hence, selecting
\[
\gamma =  \left (  \frac{ \Vert \E_{0,0} (  S^2_{k,n}  )  \Vert_q }{a_{k,n}} \right ) ^{1/(2q-2)} \, , 
\]
and noticing that we can assume that $\gamma >0$ (since otherwise there is nothing to prove), 
we get 
\begin{equation} \label{borneBq>2}
B\leq   a_{k,n}^{q-1/(2q-2)}   \Vert \E_{0,0} (  S^2_{k,n}   ) \Vert^{1/{(2q-2)}}_q \, .
\end{equation}
Starting from \eqref{forC}, considering the upper bounds \eqref{borneA} and \eqref{borneBq>2} and the fact that  
$\Vert \E_{0,0} (  S^2_{k,n}   ) \Vert \leq a_{k,n}$, we get
\begin{equation} \label{borneCq>2}
C\leq   a_{k,n}^{q-1/(2q-2)}   \Vert \E_{0,0} (  S^2_{k,n}   ) \Vert^{1/{(2q-2)}}_q \, .
\end{equation}
We handle now the quantity $E$. For  $\eta >0$, write
\[
E \leq  2^{-1} \eta \E \Big (    
(  \E_{0,n} (  S^2_{k,n}   )  )   \E^{q-1}_{0,2n} (   S^2_{k,n}     )  \Big ) + 2^{-1} \eta^{-1}  G  \, .
\]
Taking into account  \eqref{borneG} and stationarity,  it follows that 
\[
E \leq  2^{-1} \eta a^q_{k,n}    + 2^{-1} \eta^{-1}  a_{k,n}^{q-1/(q-1)}   \Vert \E_{0,0} (  S^2_{k,n}   ) \Vert^{1/{(q-1)}}_q  \, .
\]
Selecting $\eta = a_{k,n}^{-1/(2q-2)}  \Vert \E_{0,0} (  S^2_{k,n}   ) \Vert^{1/{(2q-2)}}_q$, we derive that 
\begin{equation} \label{borneE}
E\leq   a_{k,n}^{q-1/(2q-2)}   \Vert \E_{0,0} (  S^2_{k,n}   ) \Vert^{1/{(2q-2)}}_q \, .
\end{equation}
It remains to give an upper bound for $F$.  As before, for  $\eta >0$, write
\[
F \leq  2^{-1} \eta \E \big [    
(  \E_{0,2n} (  {\bar S}^2_{k,n}    )  )^q     \big ] + 2^{-1} \eta^{-1}  \E \big [    
  \E_{0,2n} (  {\bar S}^2_{k,n}    ) (  \E_{0,2n} (  {S}^2_{k,n}    )  )^{q-1}      \big ]   \, .
\]
Using the fact that $   \E_{0,2n} (  {S}^2_{k,n}    )  =  \E_{0,n} (  {S}^2_{k,n}    ) $, by the property of commuting filtrations,  it follows that 
\[
\E \big [    
  \E_{0,2n} (  {\bar S}^2_{k,n}    ) (  \E_{0,2n} (  {S}^2_{k,n}    )  )^{q-1}      \big ]  =  \E \big [    
  \E_{0,n} (  {\bar S}^2_{k,n}    ) (  \E_{0,n} (  {S}^2_{k,n}    )  )^{q-1}      \big ]  \, .
\]
So, by H\"older's inequality and stationarity, it follows that 
\[
F \leq  2^{-1} \eta  a_{k,n}^q + 2^{-1} \eta^{-1} a^{q-1}_{k,n}   \Vert \E_{0,0} (  S^2_{k,n}   ) \Vert_q \, .\]
Hence selecting $\eta = a_{k,n}^{-1/2}  \Vert \E_{0,n} (  S^2_{k,n}   ) \Vert^{1/2}_q$, we derive
\begin{equation} \label{borneF}
F\leq   a_{k,n}^{q-1/2}   \Vert \E_{0,0} (  S^2_{k,n}   ) \Vert^{1/2}_q \, .
\end{equation}
Starting from \eqref{decompositionaknq>2}, considering the upper bounds \eqref{borneA},  \eqref{borneD},  \eqref{borneG},  \eqref{borneBq>2} \eqref{borneCq>2}, \eqref{borneE} and \eqref{borneF} and taking into account the fact that $q>2$ and that $\Vert \E_{0,0} (  S^2_{k,n}   ) \Vert_q
\leq a_{k,n}$, inequality \eqref {decompositionakn2q>2} follows.  This ends the proof of the lemma. $\square$

\medskip

\noindent {\bf Proof of Proposition \ref{extensionmaxdyadic}.} We do the proof in case $d=2$. The extension for $d >2$ uses the same arguments. We proceed as in the proof of Proposition 2 in \cite{MP13}. We start by writing
\begin{multline*}
S_{2^r - m , 2^s -n } = \E ( S_{2^r  , 2^s } | {\mathcal F}_{2^r - m , 2^s -n } ) - \E ( S_{2^r  , 2^s } - S_{2^r  , 2^s -n } | {\mathcal F}_{2^r - m , 2^s -n } ) \\
-  \E ( S_{2^r  , 2^s } - S_{2^r-m  , 2^s  } | {\mathcal F}_{2^r - m , 2^s -n } ) + \E \big ( R_{m,n}(r,s)  | {\mathcal F}_{2^r - m , 2^s -n } \big ) \, ,
\end{multline*}
where $R_{m,n}(r,s) =  \sum_{k=2^r-m+1}^{2^r} \sum_{\ell=2^s-n+1}^{2^s} X_{k, \ell} $. 
Since the filtrations are commuting, by Cairoli's maximal inequality, 
\begin{multline*}
\big  \Vert \max_{0 \leq m \leq 2^r -1 \atop{0 \leq n \leq 2^s -1}} |S_{2^r - m , 2^s -n } |  \big  \Vert_p \leq q^2   \Vert  S_{2^r  , 2^s } \Vert_p  + q 
 \big   \Vert \max_{0 \leq n \leq 2^s -1} |  \E ( S_{2^r  , 2^s } - S_{2^r  , 2^s -n } | {\mathcal F}_{2^r  , 2^s -n } )  \big \Vert_p \\
+ q \big  \Vert \max_{0 \leq m \leq 2^r -1 }  \E ( S_{2^r  , 2^s } - S_{2^r-m  , 2^s  } | {\mathcal F}_{2^r - m , 2^s  } )  \big  \Vert_p +  \big  \Vert \max_{0 \leq m \leq 2^r -1 \atop{0 \leq n \leq 2^s -1}}  \big | \E  (  R_{m,n}(r,s)   | {\mathcal F}_{2^r - m , 2^s -n }   )  \big | \big \Vert_p \, .
\end{multline*}
According to the proof of Proposition 3.14 in \cite{MPU19}, using stationarity and the fact that the filtrations are commuting, we get 
\[
\big  \Vert \max_{0 \leq n \leq 2^s -1} |  \E ( S_{2^r  , 2^s } - S_{2^r  , 2^s -n } | {\mathcal F}_{2^r  , 2^s -n } )  \big \Vert_p \leq q \sum_{\ell =0}^{s-1}    2^{(s-\ell)/p}  \Vert \E ( S_{2^r,2^\ell} |  {\mathcal F}_{\infty,0})  \Vert_p 
\]
and
\[
\big  \Vert \max_{0 \leq m \leq 2^r -1} | \E ( S_{2^r  , 2^s } - S_{2^r-m  , 2^s  } | {\mathcal F}_{2^r - m , 2^s  } )  \big  \Vert_p \leq q \sum_{k =0}^{r-1}  2^{(r-k)/p}   \Vert \E ( S_{2^k,2^s} |  {\mathcal F}_{0, \infty})  \Vert_p  \, .
\]
Next, setting \[
A_{k, \ell} (r,s) =   \max_{1 \leq a \leq 2^{r-k}, a \, odd \atop{1 \leq b \leq 2^{s-\ell}, b \, odd }}   \Big | \E \Big (  \sum_{i=2^r - a2^k+1}^{2^r - (a-1)2^k} \sum_{j=2^s - b 2^\ell+1}^{2^s - (b-1)2^\ell} X_{i, j} | {\mathcal F}_{2^r - a2^k,2^s - b2^\ell} \Big )  \Big |  \, ,
\]
we have 
\[
  \big | \E  (  R_{m,n}(r,s)   | {\mathcal F}_{2^r - m , 2^s -n }   )  \big |  \leq \sum_{k=0}^{r-1} \sum_{\ell=0}^{s-1} A_{k, \ell} (r,s)   \, ,
\]
implying that 
\[
 \big  \Vert \max_{0 \leq m \leq 2^r -1 \atop{0 \leq n \leq 2^s -1}}  \big | \E  (  R_{m,n}(r,s)   | {\mathcal F}_{2^r - m , 2^s -n }   )  \big | \big \Vert_p  
 \leq   \sum_{k=0}^{r-1} \sum_{\ell=0}^{s-1}   \big   \Vert \max_{0 \leq m \leq 2^r -1 \atop{0 \leq n \leq 2^s -1}}  \E ( A_{k, \ell} (r,s)   | {\mathcal F}_{2^r - m , 2^s -n }   )  \big   \Vert_p \, .
\]
(See page 84 in \cite{MPU19} for more details).  Hence, by Cairoli's maximal inequality, 
\[
 \big  \Vert \max_{0 \leq m \leq 2^r -1 \atop{0 \leq n \leq 2^s -1}}  \big | \E  (  R_{m,n}(r,s)   | {\mathcal F}_{2^r - m , 2^s -n }   )  \big | \big \Vert_p   \leq  q^2 
   \sum_{k=0}^{r-1} \sum_{\ell=0}^{s-1}  \Vert  A_{k, \ell} (r,s) \Vert_p \, .
\]
Since 
\[
A_{k, \ell} (r,s) \leq \Big ( \sum_{a=1}^{2^{r-k} -1}  \sum_{b=1}^{2^{s-\ell} -1}  \Big |  \sum_{i= a2^k +1}^{(a+1)2^k}  \sum_{j= b2^\ell +1}^{(b+1)2^\ell} \E ( X_{i,j}   | {\mathcal F}_{a2^k,b2^\ell})  \Big | \Big )^{1/p} \, , 
\]
we get by stationarity that 
\begin{multline*}
 \big  \Vert \max_{0 \leq m \leq 2^r -1 \atop{0 \leq n \leq 2^s -1}}  \big | \E  (  R_{m,n}(r,s)   | {\mathcal F}_{2^r - m , 2^s -n }   )  \big | \big \Vert_p   \leq  q^2 
   \sum_{k=0}^{r-1} \sum_{\ell=0}^{s-1}  2^{(r-k)/p}  2^{(s-\ell)/p}  \Vert \E ( S_{2^k,2^\ell} |  {\mathcal F}_{0,0})  \Vert_p  \, .
\end{multline*}
So overall, for any positive integers $r$ and $s$, 
\begin{multline} \label{prop18d=2}
\big  \Vert \max_{1 \leq m \leq 2^r  \atop{1 \leq n \leq 2^s }} |S_{m ,n } |  \big \Vert_p \leq q^2  \Vert  S_{2^r  , 2^s } \Vert_p  +q^2 \sum_{\ell =0}^{s-1}    2^{(s-\ell)/p}  \Vert \E ( S_{2^r,2^\ell} |  {\mathcal F}_{\infty,0})  \Vert_p 
 \\
+ q^2  \sum_{k =0}^{r-1}  2^{(r-k)/p}   \Vert \E ( S_{2^k,2^s} |  {\mathcal F}_{0, \infty})  \Vert_p  
 + q^2   \sum_{k=0}^{r-1} \sum_{\ell=0}^{s-1}  2^{(r-k)/p}  2^{(s-\ell)/p}  \Vert \E ( S_{2^k,2^\ell} |  {\mathcal F}_{0,0})  \Vert_p \, ,
\end{multline}
ending the proof of the proposition when $d=2$.  The extension to the case $d>2$ follows the same strategy and is left to the reader. $\square$

\medskip

\noindent {\bf Proof of Corollary \ref{extensionmax}.} Let prove it in case  $d=2$. Taking into account Lemma 3.24 in \cite{MPU19}, and the subadditivity of the two sequences $( \Vert \E ( S_{n,m} |  {\mathcal F}_{0, \infty})  \Vert_p )_{n  \geq 1}$ and 
$( \Vert \E ( S_{n,m} |  {\mathcal F}_{\infty, m})  \Vert_p )_{m  \geq 1}$ and the fact that the filtrations are commuting, we get that there exist positive constants $c_1$, $c_2$ and $c_3$ depending on $p$ such that 
\begin{multline} \label{ineincased=2}
\big  \Vert \max_{1 \leq k \leq n  \atop{1 \leq \ell \leq m }} |S_{k ,\ell } |   \big \Vert_p \leq (2q)^2  \max_{1 \leq k \leq n  \atop{1 \leq \ell \leq m }}  \Vert  S_{k  , \ell } \Vert_p  + c_1q^2  m^{1/p}  \sum_{\ell =1}^{m}     \frac{  \max_{ 1 \leq k \leq n} \Vert \E ( S_{k,\ell} |  {\mathcal F}_{\infty,0})  \Vert_p }{\ell^{1+1/p} }
 \\
+ c_2 q^2 n^{1/p } \sum_{k=1}^{n}     \frac{  \max_{ 1 \leq \ell \leq m} \Vert \E ( S_{k,\ell} |  {\mathcal F}_{\infty,0})  \Vert_p }{k^{1+1/p} }
 + c_3 q^2  n^{1/p } m^{1/p }  \sum_{k=1}^{n}  \sum_{\ell =1}^{m}     \frac{ \Vert \E ( S_{k,\ell} |  {\mathcal F}_{0,0})  \Vert_p }{k^{1+1/p}  \ell^{1+1/p}  } \, ,
\end{multline}
which is the desired inequality in case $d=2$.   $\square$

\medskip

We also need the two following lemmas \ref{boundmomentp} and \ref{conditional-lemma}. 


\begin{lemma} \label{boundmomentp} We have 
\[
 \Vert \E_{{\bf{i}} \wedge {\bf{j}}} (X'_{{\bf{i}}} )   \E_{{\bf{i}} \wedge {\bf{j}}}  (X'_{{\bf{j}}} )   \Vert^{3/2}_{3/2}  \leq  8 M \int_0^{G(A(M,{\bf{i}},{\bf{j}} ))} Q^2(u) du \, , 
\]
where $X'_{{\bf{i}}} = \varphi_M(X_{{\bf{i}}} )  - \E ( \varphi_M(X_{{\bf{i}}} ) ) $ and $A(M,{\bf{i}},{\bf{j}} )= (4M)^{-1} \Vert \E_{{\bf{i}} \wedge {\bf{j}}} (X'_{{\bf{i}}} )   \E_{{\bf{i}} \wedge {\bf{j}}}  (X'_{{\bf{j}}} )   \Vert_1$.  
In addition
\[
\Vert \E_{{\bf{0}}} ( \varphi_M(X_{{\bf{i}}} ) \varphi_M (X_{{\bf{j}}} ) ) -  \E( \varphi_M(X_{{\bf{i}}} ) \varphi_M (X_{{\bf{j}}} ) )  \Vert^{3/2}_{3/2}
\leq  2^{3/2}  M  \int_0^{G \big ( 2^{-1} B(M,{\bf{i}},{\bf{j}} )  \big ) }  Q^2(u)  \, , 
\]
where $B(M,{\bf{i}},{\bf{j}} )= M^{-1} \Vert \E_{{\bf{0}}} ( \varphi_M(X_{{\bf{i}}} ) \varphi_M (X_{{\bf{j}}} ) ) -  \E( \varphi_M(X_{{\bf{i}}} ) \varphi_M (X_{{\bf{j}}} ) )  \Vert_1$. 
\end{lemma}
{\bf Proof of Lemma \ref{boundmomentp}.} Let $Y = \big |  \E_{{\bf{i}} \wedge {\bf{j}}} (X'_{{\bf{i}}} )   \E_{{\bf{i}} \wedge {\bf{j}}}  (X'_{{\bf{j}}} )  \big |$. By Proposition 1 in \cite{DD03}, we have 
\[
 \Vert \E_{{\bf{i}} \wedge {\bf{j}}} (X'_{{\bf{i}}} )   \E_{{\bf{i}} \wedge {\bf{j}}}  (X'_{{\bf{j}}} )   \Vert^{3/2}_{3/2}  = 
 \E ( Y^{1/2} Y)  
\leq 4 M \int_0^{(4M)^{-1}  \E ( Y)  } Q^{1/2}_{ Y } \circ G_{(4M)^{-1}Y} (u) du \, .
\]
But, using inequality (4.6) in \cite{Rio17}, we infer that 
\[
 G^{-1}_{(4M)^{-1}Y} (x)  = \int_0^x Q_{(4M)^{-1}Y} (u) du \leq   2^{-1} \int_0^x Q_{ |  \E_{{\bf{i}} \wedge {\bf{j}}} (X'_{{\bf{i}}} ) | } (u) du  \leq   \int_0^x Q(u) du   =  G^{-1} (x) \, .
\]
Therefore
\[
 \Vert \E_{{\bf{i}} \wedge {\bf{j}}} (X'_{{\bf{i}}} )   \E_{{\bf{i}} \wedge {\bf{j}}}  (X'_{{\bf{j}}} )   \Vert^{3/2}_{3/2}  
\leq 4 M \int_0^{ G ( A(M,{\bf{i}},{\bf{j}} ) } Q^{1/2}_{ Y } Q (u) du \, .
\]
Using again inequality (4.6) in \cite{Rio17}, we infer that 
\[
 \Vert \E_{{\bf{i}} \wedge {\bf{j}}} (X'_{{\bf{i}}} )   \E_{{\bf{i}} \wedge {\bf{j}}}  (X'_{{\bf{j}}} )   \Vert^{3/2}_{3/2}  
\leq 8 M \int_0^{ G ( A(M,{\bf{i}},{\bf{j}} ) }  Q^2 (u) du \, .
\]

We prove now  the second part of the lemma. Let $Z_{{\bf{0}}} = \E_{{\bf{0}}} ( \varphi_M(X_{{\bf{i}}} ) \varphi_M (X_{{\bf{j}}} ) ) -  \E( \varphi_M(X_{{\bf{i}}} ) \varphi_M (X_{{\bf{j}}} ) )$ and $Y_{{\bf{i}},{\bf{j}}} = \varphi_M(X_{{\bf{i}}} ) \varphi_M (X_{{\bf{j}}} ) $.  By Proposition 1 in \cite{DD03}, we get 
\begin{multline*}
\Vert \E_{{\bf{0}}} ( \varphi_M(X_{{\bf{i}}} ) \varphi_M (X_{{\bf{j}}} ) ) -  \E( \varphi_M(X_{{\bf{i}}} ) \varphi_M (X_{{\bf{j}}} ) )  \Vert^{3/2}_{3/2} = 
\Big |  {\rm Cov} ( |Z_{{\bf{0}}} |^{1/2} {\rm sign} ( Z_{{\bf{0}}} ) ,  Y_{{\bf{i}},{\bf{j}}}  \Big |  \\
\leq 2 M \int_0^{2^{-1} B(M,{\bf{i}},{\bf{j}} ) } Q_{|Z_{{\bf{0}}} |^{1/2} } \circ G_{M^{-1} |Y_{{\bf{i}},{\bf{j}}} |} (u) du \, .
\end{multline*}
Now note that 
\[
 G^{-1}_{M^{-1} |Y_{{\bf{i}},{\bf{j}}} |} (x)   \leq  \int_0^x Q_{|\varphi_M(X_{{\bf{0}}} )|} (u) du \leq G^{-1} (x)\, .
 \]
Hence, by a change of variables and using Lemma 2.1 in \cite{Rio17}, 
\begin{multline*}
\Vert \E_{{\bf{0}}} ( \varphi_M(X_{{\bf{i}}} ) \varphi_M (X_{{\bf{j}}} ) ) -  \E( \varphi_M(X_{{\bf{i}}} ) \varphi_M (X_{{\bf{j}}} ) )  \Vert^{3/2}_{3/2} \\
\leq 2 M \int_0^{G \big ( 2^{-1} B(M,{\bf{i}},{\bf{j}} )  \big ) } Q_{|Z_{{\bf{0}}} |^{1/2} } Q_{M^{-1} |Y_{{\bf{i}},{\bf{j}}} |} (u) du  \\
\leq 2 \int_0^{G\big ( 2^{-1} B(M,{\bf{i}},{\bf{j}} )  \big ) } Q_{|Z_{{\bf{0}}} |^{1/2} } Q_{ |Y_{{\bf{i}},{\bf{j}}} |} (u) du  \, .
\end{multline*}
Applying H\"older's inequality, it follows that 
\begin{multline*}
\Vert \E_{{\bf{0}}} ( \varphi_M(X_{{\bf{i}}} ) \varphi_M (X_{{\bf{j}}} ) ) -  \E( \varphi_M(X_{{\bf{i}}} ) \varphi_M (X_{{\bf{j}}} ) )  \Vert^{3/2}_{3/2} \\
\leq 2 \Big (  \int_0^{1} Q^3_{|Z_{{\bf{0}}} |^{1/2} }  du \Big )^{1/3}   \Big (  \int_0^{G\big ( 2^{-1} B(M,{\bf{i}},{\bf{j}} )  \big ) } Q^{3/2}_{ |Y_{{\bf{i}},{\bf{j}}} |} (u) du \Big )^{2/3}  \\
\leq  2 \Vert \E_{{\bf{0}}} ( \varphi_M(X_{{\bf{i}}} ) \varphi_M (X_{{\bf{j}}} ) ) -  \E( \varphi_M(X_{{\bf{i}}} ) \varphi_M (X_{{\bf{j}}} ) )  \Vert^{3/2}_{1/2}    \Big ( \int_0^{G \big ( 2^{-1} B(M,{\bf{i}},{\bf{j}} )  \big ) }  Q^{3/2}_{ |Y_{{\bf{i}},{\bf{j}}} |} (u) du \Big )^{2/3} \, .
\end{multline*}
Therefore, using again Lemma 2.1 in \cite{Rio17}, we derive
\begin{multline*}
\Vert \E_{{\bf{0}}} ( \varphi_M(X_{{\bf{i}}} ) \varphi_M (X_{{\bf{j}}} ) ) -  \E( \varphi_M(X_{{\bf{i}}} ) \varphi_M (X_{{\bf{j}}} ) )  \Vert^{3/2}_{3/2} \leq   2^{3/2}   \int_0^{G \big ( 2^{-1} B(M,{\bf{i}},{\bf{j}} )  \big ) }Q^{3/2}_{ |Y_{{\bf{i}},{\bf{j}}} |} (u) du  \\
\leq   2^{3/2}   \int_0^{G \big ( 2^{-1} B(M,{\bf{i}},{\bf{j}} )  \big ) } Q^3_{|\varphi_M(X_{{\bf{0}}} )|} (u)  \leq   2^{3/2}  M  \int_0^{G\big ( 2^{-1} B(M,{\bf{i}},{\bf{j}} )  \big ) }  Q^2(u)   \, ,
\end{multline*}
which is the desired inequality. $\square$

\medskip

The next lemma is given in the general setting of commuting transformations and is needed to control the dependence coefficients \eqref{plusieurs-points} in case of two commuting expanding endomorphisms. 
\begin{lemma}\label{conditional-lemma}
Let $g, \, h$ be bounded measurable functions. For every non negative integers $k,\ell, r, s, j,q$, we have , setting $u:=(K_1^{k\vee \ell-k}K_2^{r\vee s-r}g)(K_1^{k\vee \ell-\ell}K_2^{r\vee s-s}h)$ 
\begin{equation}\label{conditional-bis}
\E( U_1^kU_2^rg U_1^\ell U_2^s h |\F_{k\vee \ell+j ,r \vee s+q})
=   U_1^{k\vee \ell+j}U_2^{r\vee s+q} K_{1}^jK_{2}^qu \, .
\end{equation}
In particular,  setting  ${\tilde u}= u - \E(u)$, 
\begin{equation}\label{conditional-identity}
\Vert \E( U_1^kU_2^rg U_1^\ell U_2^s h |\F_{k\vee \ell+j ,r \vee s+q}) - \E( U_1^kU_2^rg U_1^\ell U_2^s h)  \Vert_1 
= \Vert   K_1^{j}K_2^q   \tilde u  \Vert_1 \, . 
\end{equation}
\end{lemma} 
{\bf Proof of Lemma \ref{conditional-lemma}.} The proof follows from the following considerations.  First, we notice that in addition to the property \eqref{propcom},  we also have the following properties for every 
${\mathcal A}$-measurable functions $g, \, h\, :\, \Omega\to \R$ 
such that $g,\, h\ge 0$ of $g\in {\mathbb L}^{p}(\Omega)$ and $h\in L^q(\Omega)$ with $1\le p,q\le \infty$ and $1/p+1/q=1$:
\begin{equation}\label{identity}
K_i^\ell(U_i^\ell g h)=gK_i^\ell h\qquad \forall i\in \{1,\ldots , d\},\, \ell\in \N\, .
\end{equation}

Indeed, for $f,g,h$ non-negative measurable functions we have 
\begin{gather*}
\E( fK_i^\ell(U_i^\ell g h) )= \E(U_i^\ell f U_i^\ell g h) 
=\E(U_i^\ell (fg) h)=\E(fg K_i^\ell h)\, ,
\end{gather*}
and \eqref{identity} follows for non-negative measurable functions. The general  case follows by approximating functions in ${\mathbb L}^{p}$ with simple functions.

\medskip

In particular, for every non-negative integers, we have 
\begin{equation}\label{identity-transfer}
K_i^{k\vee \ell}(U_i^k gU_i^\ell h)= K_i^{k\vee \ell-k}g
 K_i^{k\vee \ell-\ell}h\, .
\end{equation}
Indeed, using \eqref{isometry}, we see that for $k\ge \ell$,
$$
K_i^{k\vee \ell}(U_i^k gU_i^\ell g)=K_i^k U_i^\ell(U_i^{k-\ell}g \, h)=K_i^{k-\ell}(U_i^{k-\ell}g \, h),
$$
and the case where $k\le \ell$ may be proved similarly.





Finally, recall \eqref{conditional-expectation}:  for every non-negative integers $m,n$,
we have 
$$\E(\cdot |\F_{m,n})=\E\big( \E(\cdot |\F^{(2)}_n)|\F^{(1)}_m\big)= \E\big( \E(\cdot |\F^{(1)}_m)|\F^{(2)}_n\big)=
U_1^mU_2^n K_1^mK_2^n \, .
$$

Let us prove \eqref{conditional-bis}, and note that \eqref{conditional-identity} will follow from the fact that $U_1$ and $U_2$ are isometries of ${\mathbb L}^1(\mu)$. We have 
\begin{align*}
\E( U_1^kU_2^rg U_1^\ell U_2^s h |\F_{k\vee \ell+j ,r \vee s+q}) &  = U_1^{k \vee \ell +j}U_2^{r\vee s +q} K_1^{k\vee \ell +j}K_2^{r\vee s+q}\big(U_1^kU_2^rg U_1^\ell U_2^s h\big )\\
  & = U_1^{k \vee \ell +j}U_2^{r\vee s +q}K_1^jK_2^qu \, ,
\end{align*}
with $u$ as in the lemma. $\square$

\medskip

\noindent {\bf Acknowledgement.} The first author was supported by the ANR project “Rawabranch” number ANR-23-CE40-0008.

\end{document}